\documentclass[aos,MSNbibl,seceqn,dvips]{arximspdf}

\doi{10.1214/14-AOS1227} %kopijuoti is PTS
\volume{42}
\issue{4}
\pubyear{2014}
\firstpage{1361}
\lastpage{1393}

\makeatletter
\newtheorem{theorem}{Theorem}
\newtheorem{lemma}{Lemma}
\newproclaim{remark}{Remark}
\newtheorem{proposition}{Proposition}
\newproclaim{definition}{Definition}

\makeatother

\begin{document}
\begin{frontmatter}

\title{Markov jump processes in modeling coalescent with~recombination}
\runtitle{Coalescent with recombination}

\begin{aug}
\author{\fnms{Xian} \snm{Chen}\thanksref{t1}\ead[label=e1]{chenxian@amss.ac.cn}},
\author{\fnms{Zhi-Ming} \snm{Ma}\corref{}\thanksref{t1,t2}\ead[label=e2]{mazm@amt.ac.cn}}
\and
\author{\fnms{Ying} \snm{Wang}\thanksref{t1}\ead[label=e3]{wy\_math@hotmail.com}}
\runauthor{X. Chen, Z.-M. Ma and Y. Wang}
\thankstext{t1}{Supported by the National Center for Mathematics and
Interdisciplinary Sciences (NCMIS).}
\thankstext{t2}{Supported by 973 project (2011CB808000), NSFC Creative
Research Groups (11021161) and the
Fundamental Research Funds for the Central Universities (2011JBZ019).}
\affiliation{Academy of Math and Systems Science, CAS}
\address{Academy of Math and Systems Science, CAS\\
Zhong-guan-cun East Road 55\\
Beijing 100190\\
China\\
\printead{e1}\\
\phantom{E-mail:\ }\printead*{e2}\\
\phantom{E-mail:\ }\printead*{e3}}
\end{aug}

% HISTORY:
\received{\smonth{1} \syear{2014}}
\revised{\smonth{4} \syear{2014}}

% ABSTRACT
%
\begin{abstract}
Genetic recombination is one of the most important mechanisms that can
generate and maintain diversity,
and recombination information plays an important role in population
genetic studies. However, the
phenomenon of recombination is extremely complex, and hence simulation
methods are indispensable
in the statistical inference of recombination. So far there are mainly
two classes of simulation
models practically in wide use: back-in-time models and spatially
moving models. However, the statistical
properties shared by the two classes of simulation models have not yet
been theoretically studied. Based on our joint research with CAS-MPG
Partner Institute
for Computational Biology and with Beijing Jiaotong University, in this paper
we provide for the first time a rigorous argument that the statistical
properties of the two
classes of simulation models are identical. That is, they share the
same probability distribution
on the space of ancestral
recombination graphs (ARGs). As a consequence, our study provides a
unified interpretation for the algorithms of simulating coalescent
with recombination, and will facilitate the study of statistical
inference on recombination.
\end{abstract}

% KEYWORDS
\end{frontmatter}
%
%s1 #&#
\section{Introduction}\label{sec1}
Genetic recombination is an important mechanism\break  which generates and
maintains diversity. It is one of the main sources providing new
genetic materials that allow nature selection to carry on. In various
population genetic studies, such as DNA sequencing, disease study,
population history study, etc., recombination information plays an
important role. On the other hand, recombination adds much more
complexity and makes statistical inference of some evolutionary
parameters more
difficult. In the last two decades, some simulation models generating
graphs, called ancestral recombination graphs (ARGs), based on coalescent
processes have been developed to study recombination. However, none of
the existing
simulation models is perfect and each has its own advantages and
disadvantages.

Historically, a model generating the genealogical relationship
between $k$ sampled sequences from a population with constant size
without recombination was first described by Watterson
(cf. \cite{Watterson1975}), and further developed into the theory of
the coalescent by Kingman (cf. \cite{Kingman,Kingman2}). A model
describing the evolution of infinite-site sequences subject to both
coalescence and recombination in a population was first introduced
by Hudson (cf. \cite{Hudson83}). In his setup, a combined coalescent
and recombination process is followed back in time until all
nucleotide positions in the extant sequence have one common ancestral
nucleotide. The resulting structure is no longer a tree but a graph,
which was later named as ARG by
Griffiths and Marjoram, who gave in \cite{Griffiths97} more
details on ARG and embedded ARG in a birth--death process with
exponentially distributed and independent waiting times for
coalescent and recombination events. The ARG described by Griffiths
and Marjoram is simple but in many cases unnecessarily time
consuming to simulate. Because an ``ancestral'' sequence in the
birth--death process may have no genetic material in common with a sequence
descended from it. Adjusting the above shortcoming, Hudson proposed
a more efficient algorithm $ms$ (cf. \cite{Hudson02}) which is now a
commonly used computer program to simulate coalescence. Hudson's
program generates ARG back in time from the present. Due to the
Markov property of the process, the algorithm is computationally
straightforward and simple. But it is not possible to reduce the
computation further, and it is hard to approximate the computation. On the
other hand, Wiuf and Hein proposed an alternative algorithm that
moves along the sequencs and modifies the genealogy as recombination
breakpoints are encountered (cf. \cite{Wuif99}). It begins with a
coalescent tree at the left end of the sequence, and adds more
different local trees gradually along the sequence, which form part
of the ARG. The algorithm terminates at the right end of the
sequence when the full ARG is determined. Wiuf and Hein's algorithm
will produce some redundant branches in ARG. Its performance is not
so good in comparison with $ms$. But the spatially moving program is
easier to approximate. Based on the idea of constructing ARG along
sequences, there have been some approximation algorithms, such as
SMC, SMC', and MaCS (cf. \cite{McVean05,Marjoram06,Gary09}). Wiuf and Hein's spatial approach of simulating
genealogies along a sequence has a complex non-Markovian structure
in that the distribution of the next genealogy depends not just on
the current genealogy, but also on all previous ones. Therefore, the
mathematical formulation of spatial algorithm is cumbersome, and up
to date all the comparisons and discussions between spatial
algorithms and $ms$ (back-in-time algorithm) are based on simulation
studies. There is no rigorous argument showing that the ARG generated
by a spatial algorithm can share the same probability distribution
as the ARG generated by a back-in-time algorithm.

In our recent joint research with scientists in computational biology, we
proposed a new model describing coalescent with recombination, and developed
a new algorithm based on this new model. Our algorithm is also a
spatial algorithm. But we have improved Wiuf and Hein's program in that
our algorithm does not produce any redundant branches which are inevitable
in Wiuf and Hein's algorithm. In generating ARGs, our algorithm has
comparable performance with the algorithm $ms$. In addition,
our method can generate ARGs that are consistent with the sample
directly. Moreover, we can show that the existing approximation methods
(SMC, SMC', MaCS) are all special cases of our algorithm. For details,
see our
joint paper \cite{wang}.

In this paper, we further study the statistical properties of
our new model. In particular, we prove rigorously that the
statistical properties of the ARG generated by our spatially moving
model and that generated by a back-in-time model are identical, that
is, they share the same probability distribution on the space of
ARGs (cf. Theorem~\ref{coincide} below). Since the existing
approximations by spatial methods (SMC, SMC',
MaCS) are all special cases of our algorithm, consequently our study
provides a
unified interpretation for the algorithms of simulating coalescent
with recombination, and will facilitate the study of statistical
inference of recombination.

The remainder of this paper is organized as follows. As a necessary
preparation, in Section~\ref{preparation} we investigate in detail the
back-in-time model. In Section~\ref{statespace}, we describe briefly
a typical back-in-time model for simulating coalescent processes with
recombination. Then we study the state space of the Markov jump process
behind the model. In Section~\ref{construction}, we construct a
Markov jump process corresponding to the typical back-in-time
algorithm. In Section~\ref{constructionG}, we show that with
probability one, a path of the Markov jump process constitutes an ARG.
We then explore some properties of the space $G$ of ARGs. It is worth
pointing out that although Section~\ref{preparation} is a necessary
preparation, indeed our investigation is new and the results obtained
in this section have interests by their own. In particular, we believe
that the probabilistic ARG space $(G,\mathcal{B}(G),P)$ obtained in
Theorem~\ref{space} will be very useful elsewhere. In Section~\ref{mainsection}, we present our main results on the spatially moving
model. In Section~\ref{sequence},
we define and study a random sequence
$\{(S_i, Z^i),i\geq0\}$ on the probabilistic ARG space
$(G,\mathcal{B}(G),P)$, which is important for modeling our spatial
algorithm. We first define the random sequence $\{(S_i, Z^i),i\geq0\}$
and study its measurable structure. We then discuss and derive the
distributions of $S_i$ and~$Z^i$.\vadjust{\goodbreak} The derivations of the successive
conditional distributions of the involved random variables are very
complicated. Some cumbersome derivations are moved to the supplementary
article \cite{supplement}. In Section~\ref{algorithm}, we describe our
model of spatially moving
algorithm and study its statistical property. We first briefly describe
the algorithm $\mathit{SC}$ (Sequence Coalescence simulator) of our spatial
model. Afterward, we give some explanation of the algorithm. Finally,
we reach the main goal of the paper. We show that the statistical property
of the ARG generated by our spatially moving model recently proposed in
\cite{wang} is identical with the one generated by the typical
back-in-time model discussed in Section~\ref{preparation}, that is,
they share the same probability distribution on the ARG space (see
Theorem~\ref{coincide}). In Section~\ref{proofs}, we present the proofs
of the main results along with some technical lemmas. To reduce the
length of the paper, the proofs of the results in Section~\ref{preparation}
as well as part of the results in Section~\ref{mainsection} are moved to
the supplementary article~\cite{supplement}.

%s2 #&#
\section{Preparation: Investigation on the back-in-time model}\label
{preparation}
%s2.1 #&#
\subsection{State space of back-in-time model}\label{statespace}
We start our discussion by describing the state space of the Markov jump
process behind the model of a typical back-in-time algorithm.
Following Griffiths and Marjoram (cf. \cite{Griffiths96,Griffiths97}), in our model a gene or
a DNA sequence is represented by the unit interval $[0,1)$. The
model is derived from a discrete Wright--Fisher model in which, when
looking back in time,
children in a generation choose one parent with probability $1-r$,
or two parents with probability $r$; the latter means that a
recombination event occurs. If
recombination occurs, a position $S$ for the breakpoint is chosen
(independent from other breakpoints) according to a given
distribution, and the child gene is formed with the gene segments
$[0,S)$ from the first parent and $[S,1)$ from the second parent. A
continuous time model is obtained by first fixing the population size
$2N_e$ and then letting $N_e\rightarrow\infty$. Time is measured in
units of $2N_e$ and the recombination rate per gene per generation $r$
is scaled by holding $\rho=4N_er$ fixed. The limit model is a
continuous time Markov process with state space described as below.

Let $\mathcal{P}_N$ be the collection of all the subsets of
$\{1,2,\ldots,N\}$. We endow $\mathcal{P}_N$ with the discrete metric
which will be denoted by $d_p$. Throughout this paper, we shall fix $N$ and
shall hence simply write $\mathcal{P}$ for $\mathcal{P}_N$. We
denote by $\mathcal{S}_{[0,1)}(\mathcal{P})$ the family of all the
$\mathcal{P}$-valued right continuous piecewise constant functions
on $[0,1)$ with at most finitely many discontinuity points. An
element $f \in\mathcal{S}_{[0,1)}(\mathcal{P})$ may be expressed as
$f=\sum_{i=0}^m f(a_i)I_{[a_i,a_{i+1})}$ with $0=a_0<a_1<\cdots
<a_m<a_{m+1}=1$, which means that $f$ takes value $f(a_i)\in
\mathcal{P}$ on the semiclosed interval $[a_i,a_{i+1})$ for each
$i$.
%
%de1 #&#
\begin{definition}\label{state}
A finite subset $x=\{f_1,f_2,\ldots,f_k\}$ of
$\mathcal{S}_{[0,1)}(\mathcal{P})$ is said to be a state, and is
denoted by $x\in E$, if and only if $\{f_j(s)\dvtx f_j(s)\neq
\varnothing,j=1,2,\ldots,k\}$ form a partition of $\{1,2,\ldots,N\}$ for
each $s\in[0,1)$, and $f_j\not\equiv\varnothing$ for each
$j=1,2,\ldots,k$.
\end{definition}
The totality $E$ of all the states will serve as a state space of
the Markov process behind our algorithm. The process takes values in
$E$, starts at the present and traces back in time. At the present time
$X(0)= \varpi$, here
%
%e2.1 #&#
\begin{equation}\qquad
\label{x0}\varpi:=(h_1,h_2,
\ldots,h_N) \qquad\mbox{with } h_j=\{j\}I_{[0,1)}
\mbox{ for each } j=1,2, \ldots,N,
\end{equation}
representing that
the algorithm starts from $N$ sample lineages of $\mathit{DNA}$ sequences.
Starting from the present and looking back in time, if $X(t)= x$
with $x=(f_1,f_2,\ldots,f_k)\in E$, then it represents that at time
point $t$, there are $k$ ancestral lineages (i.e., there are $k$
lineages which carry the ancestral materials of the samples); if
$f_j\in X(t)$ is expressed as $f_j=\sum_{i=0}^m
f_j(a_i)I_{[a_i,a_{i+1})}$, then, for $i=1,2,\ldots,m$, on the loci
located in the
interval $[a_i,a_{i+1})$, the $j$th lineage carries ancestral
materials of the sample sequences $f_j(a_i)$. The $k$ ancestral
lineages are kept unchanged until coalescence or recombination event
happens. When a coalescence event happens, the algorithm chooses two
lineages randomly from the $k$ lineages and merges them into one
lineage. When a recombination event happens, it draws a lineage
randomly from the $k$ lineages and splits it into two lineages with
breakpoint $s$, whereas the breakpoint $s\in(0,1)$ is chosen
according to a given distribution. The waiting times between two
events are designed to be exponentially distributed with parameters
depending on the current states. The algorithm starts at the present
and performs back in time generating successive waiting times
together with recombination or coalescence events. In Griffiths and
Marjoram's program, the algorithm will repeat the above procedure
until there is only one ancestral lineage, that is, until the GMRCA
(grand most recent common ancestor) is found. To avoid redundant
computation, $ms$ algorithm has improved the above procedure in some
aspects. In particular, $ms$ respects the following two rules.

(R1) When the $N$ samples have already found common ancestry in $[0,u)$,
the algorithm will not perform any recombination event with
breakpoint located in $[0,u)$.

(R2) A locus $u$ in the sequence of a
lineage can be chosen as a breakpoint only if both $[0,u)$ and
$[u,1)$ carry ancestral materials of the samples.

In this paper, we
shall refer the algorithm with the above two
rules as a \textit{typical back-in-time algorithm}. (Remark: A minor
difference between $ms$ and our typical back-in-time algorithm is that
in the above rule (R1), $ms$ will respect common
ancestry in $[a,u)$ for any $0\leq a < u$, while we respect only
common ancestry in $[a,u)$ with $a=0$.) Since in the infinite allele
model mutation is independent of the coalescent process with
recombination, in our model we shall temporarily not consider mutation.

Below we shall show that $E$ equipped with a suitable metric $d$ is a
locally compact separable space, which is crucial and very convenient
for our further discussion.

Let $x=\{f_1,f_2,\ldots,f_k\}\in E$. We say that $s\in(0,1)$ is a
breakpoint of $x$, if there exist at least one $f_j\in x$ such that
$f_j$ is discontinuous at $s$. Suppose that $a_1<a_2<\cdots<a_m$
are the different breakpoints of $x$. We make the convention that
$a_0=0$ and $ a_{m+1}=1$, then each $f_j\in x$ may be expressed as
$f_j=\sum_{i=0}^m f_j(a_i)I_{[a_i,a_{i+1})}$. We define
%
%e2.2 #&#
\begin{equation}
\label{d} d_0(x)= \min_{0\leq i \leq m}(a_{i+1}-a_i).
\end{equation}
For $m\geq0$, we set
%
%e2.3 #&#
\begin{equation}
\label{Vm} V_m=\{x\in E\dvtx x \mbox{ has exactly } m
\mbox{ different breakpoints}\}
\end{equation}
and
%
%e2.4 #&#
\begin{equation}
\label{emstar} V_m^*=\{x\in V_m\dvtx \mbox{all the
breakpoints of } x \mbox{ are rational numbers}\}.
\end{equation}
Let $|x|$ be the number of elements in $x$. We set
%
%e2.5 #&#
\begin{equation}
\label{Uk} U_k= \bigl\{x\in E\dvtx |x|=k\bigr\}.
\end{equation}

%pr1 #&#
\begin{proposition}\label{prop1}
\textup{(i)} For $f,h
\in\mathcal{S}_{[0,1)}(\mathcal{P})$, we define
\[
d_L(f,h)=\int_0^1
d_p\bigl(f(s),h(s)\bigr)\,ds,
\]
then $d_L$ is a metric on $\mathcal{S}_{[0,1)}(\mathcal{P})$.

\textup{(ii)} Let $x,y \in E$. Suppose that $x=\{f_1,\ldots,f_k\}\in V_m\cap U_k$
and
$y=\{h_1,\ldots,\break   h_l\}\in V_n \cap U_l$, we define
%
%e2.6 #&#
\begin{equation}
\label{dxy}  d(x,y)=\max\Bigl\{\max_{1\leq i \leq k}d_L(f_i,y),
\max_{1\leq j \leq
l}d_L(h_j,x)\Bigr\}+|k-l|+
|m-n|
\end{equation}
[$d_L(f_i,y)$ stands for the distance from $f_i$ to the set $y$], then
$d$ is a metric on~$E$.
\end{proposition}
The proof of Proposition~\ref{prop1} is presented in Section A.1 of
the supplemental article \cite{supplement}.

For our purpose, we shall sometimes put the functions of a state
$x\in E$ in a parentheses to indicate that they have been well
ranked in a specific order. More precisely, we may sometimes write
$x=(f_1,f_2,\ldots,f_k)$, which means that
$(f_1,f_2,\ldots,f_k)$ have been ranked in such a way that
if $i<j$, then either
%
%e2.7 #&#
\begin{equation}
\label{order1} \min\bigl\{s\dvtx f_i(s)\neq\varnothing\bigr\} < \min
\bigl\{s\dvtx f_j(s)\neq\varnothing\bigr\}
\end{equation}
or, in case that $\min\{s\dvtx f_i(s)\neq\varnothing\} = \min
\{s\dvtx f_j(s)\neq\varnothing\} =: S$,
%
%e2.8 #&#
\begin{equation}
\label{order2} \min f_i(S) < \min f_j(S).
\end{equation}
Note that with the above specific order, the subscript ``$j$'' of an
element $f_j \in x$ is uniquely determined.

The proposition below explores
the neighborhood of a state $x\in E$.
%
%pr2 #&#
\begin{proposition}\label{neighborhood}
Let $x=(f_1,f_2,\ldots,f_k)\in V_m\cap U_k$ and $a_1<a_2<\cdots
<a_m$ be the different breakpoints of $x$.

\textup{(i)} Suppose that $y\in E$ satisfies $d(y,x)< 1$, then $y\in
V_m\cap U_k$.

\textup{(ii)} Suppose that $y=(h_1,h_2,\ldots,h_k)\in E$ satisfies
$d(y,x) < \varepsilon\leq3^{-1}d_0(x)$. Let $b_1<b_2<\cdots<b_m$ be
the different breakpoints of $y$. Then $|b_i-a_i|<\varepsilon$ for all
$1\leq i\leq m$. Moreover, for
each $1\leq j \leq k$, it holds that $d_L(h_j, f_j)= d_L(h_j,
x)=d_L(y,f_j)$, and
$h_j \in y$ can be expressed as
%
%e2.9 #&#
\begin{equation}
\label{gi} h_j=\sum_{i=0}^m
f_j(a_i)I_{[b_i,b_{i+1})}.
\end{equation}

\textup{(iii)} For any $\alpha> 0$ with $\alpha< 3^{-1}d_0(x)$, there exists an
element $y \in V_m^*$ such that $d(y,x)< \alpha$.
\end{proposition}
The proof of Proposition~\ref{neighborhood} is presented in Section
A.2 of the supplemental article \cite{supplement}.

Employing the above proposition, we can check the following topological
properties of $E$.
%
%pr3 #&#
\begin{proposition}\label{prop2} \textup{(i)} For each $m\geq0$ and
$k\geq1$, $V_m\cap U_k$ is an isolated
subset of $E$, that is, $V_m\cap U_k$ is both open and closed in
$E$.

\textup{(ii)} $E$ is a locally compact separable metric space.
\end{proposition}
Proposition~\ref{prop2} is proved in Section A.2 of the supplemental
article \cite{supplement}.
%s2.2 #&#
\subsection{Markov jump process}
\label{construction}
In this subsection, we shall construct a Markov jump process
describing the typical back-in-time algorithm. To formulate a
rigorous mathematical model, we need to introduce some operations on
$\mathcal{S}_{[0,1)}(\mathcal{P})$ and on $E$ corresponding to the algorithm.

Let $f,h \in\mathcal{S}_{[0,1)}(\mathcal{P})$. We define
\[
(f\vee h) (s):=f(s)\cup h(s).
\]
For $u \in(0,1)$, we define
\begin{eqnarray*}
f^{(u-)}(s)&:=&\cases{f(s),&\quad $\mbox{if } s< u,$\vspace*{2pt}
\cr
\varnothing,
&\quad $\mbox{if } s\geq u,$}
\\
f^{(u+)}(s)&:=&\cases{\varnothing,&\quad$\mbox{if } s<u,$\vspace*{2pt}
\cr
f(s),
&\quad$\mbox{if } s\geq u.$}
\end{eqnarray*}
Then $f\vee h, f^{(u-)}$ and $
f^{(u+)}$ are all elements of $\mathcal{S}_{[0,1)}(\mathcal{P})$.

For a state $x=(f_1,f_2,\ldots,f_k)\in E$, we set
\begin{eqnarray*}
b_1(x) &=& \inf\bigl\{u\dvtx f_1(u)\neq\{1,2,\ldots,N\}
\bigr\},
\\
b_i(x) &=& \inf\bigl\{u\dvtx f_i(u)\neq\varnothing\bigr
\}\qquad \forall 2\leq i \leq k
\end{eqnarray*}
and
\[
e_i(x) = \inf\bigl\{u\dvtx f_i(s)= \varnothing,
\forall s\in(u,1)\bigr\}\wedge1\qquad \forall 1\leq i \leq k.
\]
We define for $1\leq i\leq k$ and $u\in(b_i(x),e_i(x))$,
%
%e2.10 #&#
\begin{equation}
\label{Riu} R_{iu}(x)=\{f_j\dvtx j\neq i\}\cup\bigl
\{f_i^{(u-)},f_i^{(u+)}\bigr\},
\end{equation}
which indicates that a recombination
event happens on the lineage $f_i$ with breakpoint $u$.
Note that the definitions of $b_i$ and $e_i$ ensure that the algorithm
respects the above mentioned rule (R2). Moreover, the definition of $b_1$,
which is different from the other $b_i$, and the ranking rule specified
by (\ref{order1}) and (\ref{order2}) ensure that the algorithm respects
the rule (R1). Further, we define for $1\leq i_1<i_2\leq k$,
%
%e2.11 #&#
\begin{equation}
\label{Cii} C_{i_1,i_2}(x)=\{f_l\dvtx l\neq
i_1,i_2\}\cup\{f_{i_1}\vee f_{i_2}
\},
\end{equation}
which denotes the coalescence of the lineages $f_{i_1}$ and
$f_{i_2}$.

We can now construct a Markov jump process $\{X(t)\}$ as a rigorous
mathematical model for our typical back-in-time algorithm. In what
follows for any metric space $E$, we shall write
$\mathcal{B}(E)$ for the Borel subsets of $E$.

We define $q(x,A)$ for $x\in E$ and $A\in\mathcal{B}(E)$ as follows:
%
%e2.12 #&#
\begin{equation}\qquad
\label{qxA} q(x,A):= \sum_{1\leq i_1
< i_2\leq|x|}1_A
\bigl(C_{i_1,i_2}(x) \bigr)+\frac{\rho}{2}\sum
_{i=1}^{|x|}\int_{b_i(x)}^{e_i(x)}p(s)I_A
\bigl(R_{is}(x)\bigr)\,ds,
\end{equation}
if $|x|\geq2$, and
\[
q(x,A):=0 \qquad\mbox{if } |x|=1,
\]
where $I_A$ is the indicator function of $A$, $p(s)$ is the
density function of a given distribution on $(0,1)$, and $\rho$ is a
positive constant corresponding to a given recombination rate.
Further, we define $q(x)$ for $x\in E$ by setting
%
%e2.13 #&#
\begin{equation}
\label{qx} q(x):=q(x,E).
\end{equation}

For the terminologies involved in the proposition below, we refer to
Definition~1.9 of \cite{Chen2004}.
%
%pr4 #&#
\begin{proposition}\label{prop3}
$(q(x),q(x,A))$ defined by (\ref{qx}) and (\ref{qxA}) is a q-pair in the
sense that for each $x\in E$, $q(x,\cdot)$ is a measure on $\mathcal
{B}(E), q(x,\{x\})=0, q(x,E)\leq q(x)$; and for each $A\in\mathcal
{B}(E), q(\cdot)$ and $q(\cdot, A)$ are $\mathcal{B}(E)$-measurable.
Moreover, $(q(x),q(x,A))$ is totally stable in the sense that
$0\leq q(x)<\infty$ for all $x\in E$, and is conservative in the sense
that $q(x)=q(x,E)$ for all $x\in E$.
\end{proposition}
The proof of Proposition~\ref{prop3} is presented in Section A.3 of
the supplemental article \cite{supplement}.

By virtue of Proposition~\ref{prop3} and making use of the theory of
q-processes we
obtain the following proposition.
%
%pr5 #&#
\begin{proposition}\label{existence}
Given any initial distribution $\mu$ on $\mathcal{B}(E)$, there exists
a \mbox{q-process} $\{X(t), t\geq0\}$ corresponding to the q-pair
$(q(x),q(x,A))$, in the sense that $\{X(t)\}$ is a time homogeneous
Markov jump process satisfying: \textup{(i)} $P\{X(0) \in A\} = \mu(A)$;
\textup{(ii)} the transition probability of its embedded Markov chain is given by
%
%e2.14 #&#
\begin{equation}
\label{Pi} \Pi(x,A)=I_{\{q(x)\neq0\}}\frac{q(x,A)}{q(x)}+ I_{\{q(x)=0\}}I_A
(x);
\end{equation}
\textup{(iii)} the waiting time of its jump given $X(t)=x$ is exponentially
distributed with parameter $q(x)$.
\end{proposition}
The proof of Proposition~\ref{existence} is presented in Section
A.4 of the supplemental article \cite{supplement}.

Define
%
%e2.15 #&#
\begin{equation}
\label{Delta} \Delta:=\{f\} \qquad\mbox{with } f=\{1,2,\ldots,N\}I_{[0,1)}.
\end{equation}
Note that $\Delta$ is the only element in $U_1:= \{x\in E\dvtx   |x|=1\}$,
and hence is the only absorbing state in $E$ satisfying $q(x,E)=0$.

%pr6 #&#
\begin{proposition}\label{uniqueness}
The transition semigroup of the q-process specified by Proposition~\ref
{existence} is unique. Moreover, the
process will almost surely arrive at the absorbing state $\Delta$ in at
most finitely many jumps.
\end{proposition}
The proof of Proposition~\ref{uniqueness} is presented in Section
A.5 of the supplemental article \cite{supplement}.

%s2.3 #&#
\subsection{ARG space $G$}\label{constructionG}

Let $\{X(t), t\geq0\}$ be the Markov jump process constructed in
Proposition~\ref{existence} with initial distribution $\delta_{\{\varpi
\}}$,
where $\varpi$ is specified by (\ref{x0}). Assume that $\{X(t)\}$ is
defined on some probability space $(\Omega,\mathcal{F},P)$. Then for
each $\omega\in\Omega$,
$X(\cdot)(\omega)$ is an element in $\mathcal{S}_{[0,\infty)}(E)$,
where $\mathcal{S}_{[0,\infty)}(E)$ denotes the family of all the
$E$-valued right continuous piecewise constant functions
on $[0,\infty)$ with at most finitely many discontinuity points.
Note that not all elements $g\in\mathcal{S}_{[0,\infty)}(E)$ can be
regarded as an ARG generated by the back-in-time algorithm. Indeed, if
$g\in\mathcal{S}_{[0,\infty)}(E)$ represents an ARG, then
$g=\{g(t),t\geq0\}$ should satisfy the following two intuitive
requirements. (i) $g(0)=\varpi$ and if $g(t)\neq g(t_-)$,
then $g(t)$ is generated by a coalescent event or a recombination
event from the state $g(t_-)\in E$.
(ii) Along the path $\{g(t), t\geq0\}$, recombination events will not
happen more than once in any locus
$s\in(0,1)$. Below we shall prove that with probability one, $\{X(t),
t\geq0\}$
satisfies the above two requirements, and hence represents an ARG. To
state our result rigorously, we introduce some notation first.

For a state $x=(f_1,\ldots,f_k)\in E$ with $|x|\geq2$, we set
%
%e2.16 #&#
\begin{eqnarray}
\label{Ex} E_x&=&\bigl\{C_{i_1,i_2}(x)\dvtx 1\leq
i_1< i_2 \leq|x|\bigr\}
\nonumber
\\[-8pt]
\\[-8pt]
\nonumber
&&{}\cup\bigl\{R_{is}(x)
\dvtx
1\leq i\leq|x|, s\in\bigl(b_i(x),e_i(x)\bigr)\bigr\}.
\end{eqnarray}
For notational convenience, we shall also write $x=E_x$ if $x\in E$
with $|x|=1$.
We define a function $\mathcal{U}\dvtx  E\times E\rightarrow(0,1)\cup\{-1\}
$ by setting
%
%e2.17 #&#
\begin{eqnarray}
\label{Uxy} (x,y)&\rightarrow&\mathcal{U}(x,y)
\nonumber
\\[-8pt]
\\[-8pt]
\nonumber
&=&\cases{u,&\quad$\mbox{if there
exists a unique } u\in(0,1)$\vspace*{2pt}
\cr
&\quad$\mbox{such that }
y=R_{iu}(x) \mbox{ for some } 1\leq i \leq k;$ \vspace*{2pt}
\cr
-1,& \quad $
\mbox{else}.$}
\end{eqnarray}
We set for $g\in\mathcal{S}_{[0,\infty)}(E)$,
%
%e2.18 #&#
\begin{equation}\qquad
\label{taok} \tau_0\equiv0,\qquad \tau_n:=
\tau_n(g)=\inf\bigl\{t>\tau_{n-1}\dvtx g(t)\neq g(\tau
_{n-1})\bigr\} \qquad\forall n\geq1,
\end{equation}
with the convention that $\inf\varnothing= \infty$ and $g(\infty)=
\Delta$. For the sake of convenience, we shall write
$\mathcal{U}_n(g):=\mathcal{U}(g(\tau_{n-1}), g(\tau_n))$. We define
%
%e2.19 #&#
\begin{equation}\qquad
\label{G^} G^{\prime }:=\bigl\{g\in \mathcal{S}_{[0,\infty)}(E)\dvtx
g(\tau_0)=\varpi \mbox{ and } g(\tau_n)\in
E_{g(\tau_{n-1})} \mbox{ for all } n\geq1 \bigr\}
\end{equation}
and
%
%e2.20 #&#
\begin{equation}\qquad
\label{G} G:= \bigl\{g\in G^{\prime }\dvtx \mathcal{U}_n(g)
\neq\mathcal{U}_j(g) \mbox{ for all } n\neq j \mbox{ whenever }
\mathcal{U}_j(g)\in(0,1) \bigr\}.
\end{equation}
It is easy to see that if $g\in G$, then $g$ satisfies the above two
intuitive requirements. We shall call $G$ the
\textit{ARG space}.

%pr7 #&#
\begin{proposition}\label{thG}
There exists $\Omega_0\in\mathcal{F}$ with $P(\Omega_0)=1$, such that
for all $\omega\in\Omega_0$, we have $X(\cdot)(\omega)\in G$.
\end{proposition}
The proof of
Proposition~\ref{thG} is presented in Section A.6 of the
supplemental article \cite{supplement}.

Note that the ARG space $G$ specified by (\ref{G}) is a
subset of the $E$-valued Skorohod space $D_E[0,\infty)$. We are
going to show that $G$ equipped with the Skorohod topology is a locally
compact separable metric space.

We first introduce some terminologies and notation. For $g\in G$,
we set $\gamma(g)=\inf\{n\dvtx \tau_{n+1}(g)=\infty\}$ where $\tau_n(g)$
is defined by (\ref{taok}). Let
$Bp(g):=Bp(g(\tau_0),g(\tau_1),\ldots,g(\tau_{\gamma(g)}))$ be the
collection of all the breakpoints on $g$. Then $Bp(g)$ consists of
at most finitely many points of $(0,1)$. Moreover, by (\ref{G}) the
points of $Bp(g)$ are all different from each other. Denote by
$|Bp(g)|$ the number of points contained in $Bp(g)$. We define
$S_i:=S_{i}(g)$ to be the $i$th order statistic of $Bp(g)$. That
is, $Bp(g)=\{S_1,S_2,\ldots,S_{|Bp(g)|}\}$ and $S_i<S_{i+1}$ for all
$i$. For convenience, we make the convention that $S_0=0$ and $S_i=
1$ for $i>|Bp(g)|$. Suppose that $|Bp(g)|= m$ and $g(t)=
\{f_1,f_2,\ldots,f_{|g(t)|})\}\in E$, we define an
$(m+1)$-dimensional $\mathcal{P}$-valued vector for each $f_j\in
g(t)$ by setting
$\mathfrak{S}_j(g(t)):=(f_j(S_0),f_j(S_1),\ldots,f_j(S_m))$. Further,
we write $\mathfrak{S}(g(t))=\{\mathfrak{S}_j(g(t))\dvtx 1\leq j\leq
|g(t)|\}$. It is clear that
$\mathfrak{S}(g(t))=\mathfrak{S}(g(\tau_n))$ when $\tau_n \leq t <
\tau_{n+1}$ for all $n$. In what follows, we set $d_0(g)=\min_{0\leq
i\leq|Bp(g)|}(S_{i+1}(g)-S_i(g))$.

For the convenience of the reader, we recall the definition of the
Skorohod metric $d_S$ on $D_E[0,\infty)$ (cf. \cite{Kurtz}).

Let $\Lambda$ be the collection of Lipschitz continuous and strictly
increasing functions $\lambda$ such that $\kappa(\lambda):=\sup_{s>t\geq0}|\log\frac{\lambda(s)-\lambda(t)}{s-t}|<\infty$. For
$g_1,g_2\in\break  D_E[0,\infty)$, the Skorohod metric $d_S$ is defined as
%
%e2.21 #&#
\begin{equation}
\label{Gmetric} d_S(g_1,g_2)=\inf
_{\lambda\in\Lambda}\biggl[\kappa(\lambda)\vee\int_0^{\infty
}e^{-u}d(g_1,g_2,
\lambda,u)\,du\biggr],
\end{equation}
where
\[
d(g_1,g_2,\lambda,u)=\sup_{t\geq0}d
\bigl(g_1(t\wedge u),g_2\bigl(\lambda (t)\wedge u\bigr)
\bigr)\wedge1.
\]

The proposition below plays an important role in our further study.
%
%pr8 #&#
\begin{proposition}\label{2Gneighborhood}
Let $g_l,g_0\in G$. Suppose that
%
%e2.22 #&#
\begin{equation}
\label{l0} d_S(g_l,g_0)<
3^{-1}d_0(g_0)e^{-2\tau_{\gamma(g_0)}(g_0)}.
\end{equation}
Then the following assertions hold:
\begin{longlist}[(iii)]
\item[(i)] $\gamma(g_l)=\gamma(g_0)$.

\item[(ii)] $|Bp(g_l)|=|Bp(g_0)|$ and $\mathfrak{S}(g_l(\tau_n))=\mathfrak
{S}(g_0(\tau_n))$ for all $1\leq n\leq\gamma(g_0)$.

\item[(iii)] $d(g_l(\tau_n),g_0(\tau_n))\leq e^{2\tau_{\gamma
(g_0)}(g_0)}d_S(g_l,g_0)<3^{-1}d_0(g_0)$ for all $1\leq n\leq\gamma(g_0)$.
\end{longlist}
\end{proposition}
Proposition~\ref{2Gneighborhood} is proved in Section A.7 of the
supplemental article \cite{supplement}.

%pr9 #&#
\begin{proposition}\label{3Gneighborhood} \textup{(i)} Let $\{g_l, l\geq1\}
\subset G$ and $g_0\in G$. Then
$\lim_{l\rightarrow\infty}d_S(g_l,\break
g_0)=0$ if and only if $S_i(g_l)\rightarrow S_i(g_0)$ for all $i\geq1,
\tau_n(g_l)\rightarrow\tau_n(g_0)$ for all $n\geq1$, and there exists
$l_0$ such that for all $l\geq l_0$ the assertions of Proposition~\ref
{2Gneighborhood}\textup{(i)--(iii)} hold.

\textup{(ii)} $G$ equipped with the Skorohod metric $d_S$ is a locally compact
separable metric space.
\end{proposition}
The proof of Proposition~\ref{3Gneighborhood} is presented in
Section A.8 of the supplemental article \cite{supplement}.

Note that $G$ can be regarded as the collection of all the ARGs
generated by the back-in-time algorithm.
We denote by $\mathcal{B}(G)$ the Borel sets of $G$.
%
%th1 #&#
\begin{theorem}\label{space}
Let $P$ be the probability distribution on $(G,\mathcal{B}(G))$
generated by the typical back-in-time algorithm,
and denote by $\{X(t),t\geq0\}$ the coordinate process on $G$. Then
$\{X(t), t\geq0\}$ is an $E$-valued Markov jump process
corresponding to the q-pair (\ref{qxA})--(\ref{qx}).
\end{theorem}
\begin{pf}
Since on the Skorohod space the Borel $\sigma$-field coincides with the
$\sigma$-field generated by its coordinate process (cf., e.g., \cite
{Kurtz}), the theorem follows directly from Propositions \ref{thG}
and~\ref{3Gneighborhood}(ii).
\end{pf}
Before concluding this subsection, we explore some properties of
$Bp(g)$ and $\mathfrak{S}(g(t))$ as stated in Proposition~\ref{Sgt}
below, which will play an important role in our further discussion.

Below we denote by $\mathcal{P}^{m+1}$ the totality of
$(m+1)$-dimensional $\mathcal{P}$-valued vectors. For
$\vec{z}=(z_0,z_1,\ldots,z_m)\in\mathcal{P}^{m+1}$, we define
$\pi_j(\vec{z})=z_j$ for $0\leq j\leq m$. For $\vec{a}= (a_0,a_1,\ldots,a_m)\in
\mathcal{P}^{m+1}$ and $\vec{b}= (b_0,b_1,\ldots,b_m)\in
\mathcal{P}^{m+1}$, we define $\vec{a}\vee\vec{b}\in
\mathcal{P}^{m+1}$ by setting $\pi_j(\vec{a}\vee\vec{b})= a_j\cup
b_j$. Further, for $1\leq j \leq m$, we define $(\vec{a})^{j-}\in
\mathcal{P}^{m+1}$ by setting $\pi_i((\vec{a})^{j-})=a_i$ for $i<j$
and $\pi_i((\vec{a})^{j-})=\varnothing$ for $i\geq j$, define
$(\vec{a})^{j+}\in\mathcal{P}^{m+1}$ by setting
$\pi_i((\vec{a})^{j+})=\varnothing$ for $i<j$ and
$\pi_i((\vec{a})^{j+})=a_i$ for $i\geq j$. We say that a vector
$\vec{z} \in\mathcal{P}^{m+1}$ is null, if
$\pi_j(\vec{z})=\varnothing$ for all $0\leq j \leq m$.

%pr10 #&#
\begin{proposition}\label{Sgt}
For $g\in G$ with $|Bp(g)|= m$ and $\gamma(g)=\gamma$, we denote by
$S_i=S_{i}(g)$ the $i$th order statistic of $Bp(g)$ for $1\leq i \leq
m$, and write $\mathfrak{S}(t)$ for $\mathfrak{S}(g(t))$. Then the
following assertions hold:

\textup{(i)} For all $t$, $\mathfrak{S}(t)$ is a finite subset of $\mathcal
{P}^{m+1}$ such that $\{\pi_i(\vec{z})\dvtx \vec{z} \in\mathfrak{S}(t),\break  \pi
_i(\vec{z})\neq\varnothing\}$ form a partition of $\{1,2,\ldots,N\}$ for
each $0\leq i\leq m$. Moreover, any $\vec{z}\in\mathfrak{S}(t)$ is not null.

\textup{(ii)} There exist $\{\tau_n\dvtx 1\leq n \leq\gamma\}$ with $\tau_0:=0 <
\tau_1<\tau_2<\cdots<\tau_{\gamma}<\infty:=\tau_{\gamma+1}$ such that
$\mathfrak{S}(\tau_n)\neq\mathfrak{S}(\tau_{n+1})$ and
$\mathfrak{S}(t)=\mathfrak{S}(\tau_n)$ when $t\in[\tau_n, \tau_{n+1})$
for all $0\leq n\leq\gamma$.

\textup{(iii)} For $1\leq n \leq\gamma$, if we write $\mathfrak{S}(\tau
_{n-1})= \{\vec{z}_1,\vec{z}_2, \ldots,\vec{z}_k\}$. Then either
$ \mathfrak{S}(\tau_{n})=\{\vec{z}_l\dvtx l\neq j_1,j_2\}\cup\{\vec
{z}_{j_1}\vee\vec{z}_{j_2}\}$ for some $1\leq j_1<j_2\leq k$, or $
\mathfrak{S}(\tau_{n})=\{\vec{z}_l\dvtx l\neq j\}\cup\{(\vec
{z}_j)^{i-},(\vec{z}_j)^{i+}\}$ for some $1\leq j \leq k$ and some
$1\leq i \leq m$.

\textup{(iv)} For each $1\leq i \leq m$, there exists
$\tau_{n(i)}\in\{\tau_1,\tau_2,\ldots,\tau_{\gamma}\}$ at which $S_i$
appears in the following sense:
if we write $\mathfrak{S}(\tau_{n(i)-1})= \{\vec{z}_1,\vec{z}_2, \ldots,\vec{z}_k\}$, then
$\mathfrak{S}(\tau_{n(i)})=\{\vec{z}_l\dvtx l\neq j\}\cup\{(\vec
{z}_j)^{i-},(\vec{z}_j)^{i+}\}$ for
some $\vec{z}_j \in\mathfrak{S}(\tau_{n(i)-1})$ satisfying $\pi
_{i-1}(\vec{z}_j)= \pi_{i}
(\vec{z}_j)\neq\varnothing$. Moreover, the time point
$\tau_{n(i)}$ at which $S_i$ appears is unique.

\textup{(v)} For $t\in[\tau_n, \tau_{n+1})$, if we write $\mathfrak{S}(\tau
_{n})= \{\vec{z}_1,\vec{z}_2, \ldots,\vec{z}_k\}$ and make the
convention that $S_0=0$ and $S_{m+1}= 1$, then $g(t)$ is expressed as
$g(t)= \{f_1,f_2,\ldots,f_{k}\}$ with $f_j=\sum_{l=0}^{m}\pi_l(\vec
{z}_j)I_{[S_l,S_{l+1})}$ for each $1\leq j \leq k$.
\end{proposition}
\begin{pf}
All the assertions can be checked directly by the definition of $Bp(g)$
and $\mathfrak{S}(g(t))$, as well as Definition~\ref{state}, (\ref
{G^}) and (\ref{G}), we leave the details to the reader.
\end{pf}

%s3 #&#
\section{Main results: Spatially moving model}\label{mainsection}
%s3.1 #&#
\subsection{Random sequence \texorpdfstring{$\{(S_i,Z^i),i\ge0\}$}{\{(Si,Zi), i>=0\}}}\label{sequence}
Let $(G,\mathcal{B}(G),P)$ be the probability space specified in
Theorem~\ref{space}. In this subsection, we shall define a sequence of random
variables $\{(S_i, Z^i),i\geq0\}$ on
$(G,\mathcal{B}(G),P)$ and derive their distributions, which will be
used to model our spatial algorithm.

%s3.1.1 #&#
\subsubsection{Definition and structure of \texorpdfstring{$\{(S_i,Z^i),i\geq0\}$}{\{(Si,Zi), i>=0\}}}\label{SiZi}
In Section~\ref{constructionG}, we have defined $S_i:=S_{i}(g)$ to
be the $i$th order statistic of $Bp(g)$. With the convention that
$S_0=0$ and $S_i= 1$ for $i>|Bp(g)|$, by Proposition~\ref
{3Gneighborhood} we see that $\{S_i,i\geq0\}$ is a sequence of
continuous mappings from $G$ to $[0,1]$.
In what follows, we define the random sequence $\{ Z^i,i\geq0\}$.

Below for
$\vec{z}=(z_0,z_1,\ldots,z_m)\in\mathcal{P}^{m+1}$ and $0\leq i < m$,
we write $\pi_{[0,i]}(\vec{z})=(z_0, \ldots,z_i)$. For a subset
$A\subset\mathcal{P}^{m+1}$, we write $\pi^*_{[0,i]}(A):=
\{\pi_{[0,i]}(\vec{a})\dvtx   \vec{a} \in A,\break    \pi_{[0,i]}(\vec{a})
  \mbox{ is not null in }  \mathcal{P}^{i+1}\}$ for $i < m$, and make
the convention that\break
$\pi^*_{[0,i]}(A)= A$ for $i \geq m$.
Let $g\in G$. For $i= 0$, we define
%
%e3.1 #&#
\begin{equation}
\label{Z^O} T^0_0=0,\qquad Z^0(t):=
Z^0(t) (g)= \pi^*_{[0,0]}\bigl(\mathfrak{S}\bigl(g(t)\bigr)
\bigr).
\end{equation}
For $1\leq i\leq|Bp(g)|$, we set
%
%e3.2 #&#
\begin{equation}
\label{xii}\xi^i:=\xi^i(g)= \pi_{i}(
\vec{z}_j),
\end{equation}
where $\pi_{i}(\vec{z}_j)$ is specified in Proposition~\ref{Sgt}(iv),
that is, $\xi^i$ is the type set involved at
the recombination at locus $S_i$ when it first becomes a breakpoint.
Further, we define
%
%e3.3 #&#
\begin{equation}\qquad
\label{Z^i} T^i_0:= T^i_0(g)=
\tau_{n(i)},\qquad Z^i(t):=Z^i(t) (g)=
\pi_{[0,i]}(\vec {z}_l) \qquad\mbox{for } t\geq T^i_0,
\end{equation}
here $\tau_{n(i)}$ is specified in Proposition~\ref{Sgt}(iv) and $\vec
{z}_l$ is the unique vector in $\mathfrak{S}(g(t))$ satisfying $\pi
_i(\vec{z}_l)\supseteq\xi^i$. Note that the existence and uniqueness
of $\vec{z}_l$ employed above is ensured by Proposition~\ref{Sgt}(i).
Intuitively, $Z^i(t)$ traces those lineages containing the
genotypes $\xi^i$ at locus $S_i$, the ancestral materials $\vec{z}_l$ and
$\vec{z}_j$ are objects at different times $t$ and $\tau_{n(i)-1}$,
respectively. For $i> |Bp(g)|$, we make the convention that
$T^i_0=\infty$ and $Z^i(t)\equiv\vec{\varnothing}$, here $\vec{\varnothing
}$ denotes the null vector in $\mathcal{P}^{i+1}$.

For each $i\geq1$, we define recursively for $n\geq1$,
%
%e3.4 #&#
\begin{equation}
\label{Tin} T^i_n= \inf\bigl\{t>T^i_{n-1}
\dvtx Z^i(t)\neq Z^i\bigl({T^i_{n-1}}
\bigr)\bigr\}
\end{equation}
and
%
%e3.5 #&#
\begin{equation}
\label{Xiin} \xi^i_n=Z^i
\bigl(T^i_n\bigr),
\end{equation}
with the convention that $Z^i(\infty)=\vec{\varnothing}$. For
convenience, we make the further convention that $Z^i(t)=\vec{\varnothing
}$ when $t< T_0^i$. Then $\{Z^i(t), t\geq0\}$ is uniquely determined
by the $(\overline{R}^{+}\times\mathcal{P}^{i+1})$-valued sequence $\{
(T^i_n,\xi^i_n )\dvtx  n\geq0\}$. (Here and henceforth, $\overline
{R}^{+}:=[0, \infty]$.) We remind the reader that for $1\leq i\leq
|Bp(g)|$, we have $\xi^i_0=(\varnothing,\ldots,\varnothing,\xi^i)$ where
$\xi^i$ was used as a label for defining $Z^i(t)$ [cf. (\ref{xii}), (\ref
{Z^i})]. Below we endow the product topology on $\overline{R}^{+}\times
\mathcal{P}^{i+1}$.
%
%pr11 #&#
\begin{proposition}\label{Zjcon}\textup{(i)}
For each $i$ and $n$, $(T^i_n,\xi^i_n)$ is a continuous functional from
$G$ to $\overline{R}^{+}\times\mathcal{P}^{i+1}$.

\textup{(ii)} For each $i$, $\{Z^i(t), t\geq0\}$ is a jump process on $G$ with
at most finitely many jumps.
\end{proposition}
Proposition~\ref{Zjcon} is proved in Section A.9 of the supplemental
article \cite{supplement}.

We now study the measurable structure of $\{(S_i, Z^i),i\geq0\}$.
Let $g\in G$. We define for $i=0$,
%
%e3.6 #&#
\begin{equation}
\label{VZ^0} V\bigl(Z^0;t\bigr)=Z^0(t) (g).
\end{equation}
For $1 \leq i \leq Bp(g)$, we define $V(Z^0,Z^1,\ldots,Z^i;t):=
V(Z^0,Z^1,\ldots,Z^i;t)(g)$ recursively by the scheme below.
For $t< T^i_0$, define
%
%e3.7 #&#
\begin{eqnarray}\qquad
\label{t<T} &&V\bigl(Z^0,Z^1,\ldots,Z^i;t
\bigr)
\nonumber
\\[-8pt]
\\[-8pt]
\nonumber
&&\qquad=
\bigl\{\vec{z}\in\mathcal{P}^{i+1}\dvtx \pi_{[0,i-1]}(\vec{z})
\in V\bigl(Z^0,Z^1,\ldots,Z^{i-1};t\bigr),
\pi_i(\vec{z})= \pi_{i-1}(\vec{z}) \bigr\};
\end{eqnarray}
for $t\geq T^i_0$, define
%
%e3.8 #&#
\begin{eqnarray}
\label{t>T}&& V\bigl(Z^0,Z^1,\ldots,Z^i;t
\bigr)\nonumber\\
&&\qquad=\bigl\{Z^i(t)\bigr\}
\nonumber
\\[-8pt]
\\[-8pt]
\nonumber
&&\qquad\quad{} \cup
\bigl\{\vec{z}\in\mathcal{P}^{i+1}\dvtx \pi_{[0,i-1]}(\vec{z})
\in V\bigl(Z^0,Z^1,\ldots,Z^{i-1};t\bigr)
\setminus\bigl\{\pi_{[0,i-1]}\bigl(Z^i(t)\bigr)\bigr\},\\
&&\hspace*{220pt}\qquad\quad\pi_i(\vec{z})= \pi_{i-1}(\vec{z})\setminus
\xi^i \bigr\}.
\nonumber
\end{eqnarray}
For $i> Bp(g)$, we define $V(Z^0,Z^1,\ldots,Z^i;t)=\mathfrak{S}(g(t))$.
%
%pr12 #&#
\begin{proposition}\label{V}
Let $V(Z^0,Z^1,\ldots,Z^i;t):= V(Z^0,Z^1,\ldots,Z^i;t)(g)$ be defined
as above. Then for each $i \geq0$, we have
%
%e3.9 #&#
\begin{equation}
\label{VS} V\bigl(Z^0,Z^1,\ldots,Z^i;t
\bigr)=\pi^*_{[0,i]}\bigl(\mathfrak{S}\bigl(g(t)\bigr)\bigr).
\end{equation}
\end{proposition}
Proposition~\ref{V} is proved in Section A.10 of the supplemental
article \cite{supplement}.

Next, for $s\in[0,1)$ and $f \in\mathcal{S}_{[0,1)}(\mathcal{P})$, we
define $f^s \in\mathcal{S}_{[0,1)}(\mathcal{P})$ by setting
\[
f^{s}(u):= \cases{f(u), &\quad$\mbox{if } u< s,$\vspace*{2pt}
\cr
f(s),&\quad
$\mbox{if } u\geq s.$}
\]
For $x=\{f_1,f_2,\ldots,f_k\}\in E$, we define $\pi_{[0,s]}^{E}(x) \in
E$ by setting
%
%e3.10 #&#
\begin{equation}
\label{piE} \pi_{[0,s]}^{E}(x):=\bigl\{f^s_j
\dvtx 1\leq j \leq k, f^s_j\not\equiv\varnothing\bigr\}.
\end{equation}
Applying Proposition~\ref{neighborhood}(ii), one can check that $\pi
_{[0,s]}^{E}$ is a measurable map from $(E, \mathcal{B}(E))$ to $(E,
\mathcal{B}(E))$.
Below we shall sometimes write $x^s=\pi_{[0,s]}^{E}(x)$, and write
$[x^s]= (\pi_{[0,s]}^{E})^{-1}\{x^s\}=:\{y\in E\dvtx y^s = x^s\}$. Let
$\sigma(\pi_{[0,s]}^{E})$ be the sub $\sigma$-algebra of $\mathcal
{B}(E)$ generated by $\pi_{[0,s]}^{E}$. Then $[x^s]$ is an atom of
$\sigma(\pi_{[0,s]}^{E})$ for each $x\in E$.

For $g\in G$, we define $\pi_{[0,s]}^{G}(g)$ by setting
%
%e3.11 #&#
\begin{equation}
\pi_{[0,s]}^{G}(g) (t):=\pi_{[0,s]}^{E}
\bigl(g(t)\bigr)\qquad \forall t\geq0.
\end{equation}

%pr13 #&#
\begin{proposition}\label{measG}
$\pi_{[0,s]}^{G}$ is a measurable map from $(G,\mathcal{B}(G))$ to
$(G,\mathcal{B}(G))$.
\end{proposition}
Proposition~\ref{measG} is proved in Section A.11 of the
supplemental article \cite{supplement}.

We extend the definition of $\pi_{[0,s]}^{G}(g)$ by setting $\pi
_{[0,s]}^{G}(g)=g$ for $s\geq1$. Write $X^{s}(g):=\pi_{[0,s]}^{G}(g)$.
Then $\{X^s, s\geq0\}$ can be viewed as a $G$-valued stochastic
process defined on the probability space $(G, \mathcal{B}(G),P)$. From
Proposition~\ref{Sgt}(v), we see that $\{X^{s}\}$ is a jump process,
that is, its pathes are piecewise constant and right continuous with
left limits. Define $S_0^{\prime }=0$ and $S_i^{\prime }=\inf\{t>S_{i-1}^{\prime }\dvtx
X^{s}\neq X^{S_{i-1}^{\prime }}\}$ for $i\geq1$. That is, $S_i^{\prime }$ is the
$i$th jump time of $\{X^s\}$.
Let $\mathcal{F}^s=\sigma(X^{u}, u\leq s), s\geq0$, be the natural
filtration of $\{X^s\}$ and $\mathcal{F}^{\infty}=\bigvee_{s\geq0}\mathcal
{F}^s$. Since for $0<u\leq s$, it holds that $X^{u}=\pi
_{[0,u]}^{G}(X^{s})$, therefore, $X^{u}$ is $\sigma(X^{s})$ measurable.
Thus, $\mathcal{F}^s=\sigma(X^{s})$ and $\{X^{s}, s\geq0\}$ is a
$G$-valued Markov process with respect to its natural filtration.
The proposition below shows that $\{(S_i, Z^i),i\geq0\}$ enjoys a very
nice measurable structure.
%
%pr14 #&#
\begin{proposition}\label{fil}For $i\geq1$, we have
\[
\sigma\bigl(S_1,\ldots,S_i;Z^0,Z^1,
\ldots,Z^i\bigr)=\sigma\bigl(X^{S_i}\bigr) =\sigma
\bigl(X^0,S_1,\ldots,S_i;X^{S_1},
\ldots,X^{S_i}\bigr).
\]
\end{proposition}
Proposition~\ref{fil} is proved in Section~\ref{prooffil}.

%s3.1.2 #&#
\subsubsection{Distribution of $S_i$}\label{Si}

We write $\pi_t(g)=g(t)$ for $g\in G \subset D_E[0,\infty)$. For fixed
$s\geq0$, we write $X^s(t)(g)=\pi_t(X^s(g))$. It is easy to see that
$X^s(t)=X(t)^s:=\pi_{[0,s]}^E(X(t))$.\vspace*{1pt}
Therefore, $\{X^s(t)(\cdot), t\geq0\}$ is a jump process taking values
in $(E^s,\mathcal{B}(E^s))$.
Here and henceforth,
%
%e3.12 #&#
\begin{equation}
\label{Es}E^s:=\pi_{[0,s]}^{E}(E)=\bigl
\{x^s\dvtx x\in E\bigr\}.
\end{equation}
Note that $E^s=\{x\in E\dvtx   x=\pi_{[0,s]}^{E}(x)\}$, hence $E^s$ is a
Borel subset of $E$.

Set $\tau_0^{s}=0$ and for $n\geq1$ define
%
%e3.13 #&#
\begin{equation}
\label{taons} \tau_{n}^{s}=\inf\bigl\{t>
\tau_{n-1}^{s}\dvtx X^{s}(t)\neq X^{s}
\bigl(\tau _{n-1}^{s}\bigr)\bigr\}.
\end{equation}
Since $\mathcal{B}(G)=G\cap\mathcal{B}(D_E[0,\infty))=G\cap
\sigma\{\pi_t^{-1}(B)\dvtx B\in\mathcal{B}(E), t\geq0\}$, we have that
%
%e3.14 #&#
\begin{eqnarray}
\label{Fs}\mathcal{F}^s&=&\sigma\bigl(X^{s}\bigr)=
\sigma\bigl\{ X^{s}(t)\dvtx t\geq 0\bigr\}
\nonumber
\\[-8pt]
\\[-8pt]
\nonumber
&=&\sigma\bigl\{X^s(0),\tau_1^{s},X^s
\bigl(\tau_1^{s}\bigr),\tau_2^{s},X^s
\bigl(\tau_2^{s}\bigr), \ldots\bigr\}.
\end{eqnarray}

The proposition below is crucial for deriving the distributions of $\{
(S_i,Z^{i})\dvtx\break  i\geq0\}$. Its proof is quite long and involves a study of
projections of q-processes. To avoid digression from the main topics of
this paper, the proof will appear elsewhere.
In what follows, we
always assume that $y_0=\varpi$ where $\varpi$ was specified by~(\ref{x0}).
%
%pr15 #&#
\begin{proposition}\label{distribution2}Let $0\leq u\leq s\leq1$. For
$n\geq1$ and $k\geq0$, we define $\vartheta_{j}:= \tau_n^{s}+\tau
_{j-n}^{u}\circ\theta_{\tau_n^{s}}$ for $n\leq j\leq n+k$, where $\theta
_{\tau_n^{s}}$ is the time shift operator with respect to the $(\mathcal
{F}_t)$-stopping time $\tau_n^{s}$ ($(\mathcal{F}_t)_{t\geq0}$ refers
to the natural filtration generated by $\{X(t), t\geq0\}$). Then for
$B\in\mathcal{B}( (R^{+}\times E)^{n+k})$, we have
\begin{eqnarray*}
&&P\bigl\{\bigl(\tau_1^{s},X^s\bigl(
\tau_1^{s}\bigr),\ldots,\tau _n^{s},X^s
\bigl(\tau_n^{s}\bigr), \vartheta_{n+1},X^u(
\vartheta_{n+1}),\ldots, \vartheta_{n+k},X^u(
\vartheta_{n+k})\in B\bigr)\bigr\}\\
&&\qquad=
\int_0^{\infty}\,dt_1\cdots\int
_0^{\infty}\,dt_{n+k}\int_Eq(y_0,dy_1)
\cdots\\
&&\qquad\quad \int_E q(y_{n-1},dy_n)\int
_E q\bigl(\pi_{[0,u]}^E(y_n),dy_{n+1}
\bigr)
\\
&&\qquad\quad \int_E q(y_{n+1},dy_{n+2})\cdots\int
_E q(y_{n+k-1}, dy_{n+k})I_{B}(t_1,y_1,
\ldots,t_{n+k},y_{n+k})
\\
&&\qquad\quad{} \cdot I_{\{t_1<\cdots< t_{n+k}\}}(t_1,\ldots,t_{n+k})\prod
_{j=0}^{n-1}I_{E^{s}}(y_{j+1})\prod
_{i=0}^{k-1}I_{E^{u}}(y_{n+i+1})
\\
&&\qquad\quad{} \cdot\exp\Biggl\{-\sum_{j=0}^{n-1}q
\bigl(y_j,\bigl[y_j^s\bigr]^{c}
\bigr) (t_{j+1}-t_j) \\
&&\hspace*{59pt}{}-\sum_{i=0}^{k-1}q
\bigl(y_{n+i},\bigl[y_{n+i}^u\bigr]^{c}
\bigr) (t_{n+i+1}-t_{n+i})\Biggr\}.
\end{eqnarray*}
\end{proposition}
\begin{pf}
See \cite{XChen2013}.
\end{pf}

For $x=\{f_1,f_2,\ldots,f_k\}\in
E$, we define $\pi_s^E(x) \subset\mathcal{P}$ by setting
%
%e3.15 #&#
\begin{equation}
\label{pis} \pi_s^E(x):=\bigl\{f_j(s)
\dvtx f_j(s)\neq\varnothing,j=1,2,\ldots,k\bigr\}.
\end{equation}
Note that by Definition~\ref{state}, $\pi_s^E(x)$ is a partition of
$\{1,2,\ldots,N\}$.
For $g\in G$, $s\in[0,1)$, we write
%
%e3.16 #&#
\begin{equation}
\label{Ts} \mathcal{T}_s(t) (g)=\pi_s^E
\bigl(g(t)\bigr) \quad\mbox{and}\quad L_s(g)=\int_0^{\beta
_s}\bigl|
\mathcal{T}_{s}(t) (g)\bigr|\,dt,
\end{equation}
where $\beta_s=\inf\{t\dvtx |X^{s}(t)(g)|=1 \}$.
Intuitively, $\{\mathcal{T}_s(t), t\geq0\}$ is the coalescent tree at
site $s$,
and $L_s(g)$ is the total length
of the coalescent tree $\mathcal{T}_{s}(g)$ before~$\beta_s$.

%th2 #&#
\begin{theorem}\label{distSi}
For $i\geq0$, the distribution of $S_{i+1}$ conditioning on $\mathcal
{F}^{S_i}$ is: for $s<1$,
\[
P\bigl(S_{i+1}> s|\mathcal{F}^{S_i}\bigr)=\exp\biggl\{-\rho
L_{S_i}\bigl(X^{S_i}\bigr)\int_{S_i}^{s\vee S_i}2^{-1}p(r)\,dr
\biggr\}
\]
and
\[
P\bigl(S_{i+1}=1|\mathcal{F}^{S_i}\bigr)=\exp\biggl\{-\rho
L_{S_i}\bigl(X^{S_i}\bigr)\int_{S_i}^{1}2^{-1}p(r)\,dr
\biggr\}.
\]
\end{theorem}
Theorem~\ref{distSi} is proved in Section~\ref{proofdistSi}.

%s3.1.3 #&#
\subsubsection{Distribution of $Z^i$}\label{subZ^i}
%
%th3 #&#
\begin{theorem}\label{distZ^0}
We have $Z^0(t)=\mathcal{T}_0(t)$, and the distribution of $\mathcal
{T}_0= \{\mathcal{T}_0(t),\break  t\geq0\}$ follows that of a standard
Kingman's coalescent tree developed in \cite{Kingman}.
\end{theorem}
Theorem~\ref{distZ^0} is proved in Section~\ref{proofdistZ^0}.

Below we study the distribution of $Z^{i+1}$ conditioning on $\mathcal
{F}^{S_i} \vee\sigma(S_{i+1})$ for each $i\geq0$. Note that for
$i\geq1$, $\{Z^{i}(t), t\geq0\}$ is uniquely determined by the
$(\overline{R}^{+}\times\mathcal{P}^{i+1})$-valued sequence $\{
(T^i_n,\xi^i_n )\dvtx  n\geq0\}$. Thus, by virtue of Proposition~\ref{fil},
we need only to calculate the distribution of $\{(T^{i+1}_n,\xi^{i+1}_n
)\dvtx  n\geq0\}$ conditioning on $\sigma(X^{S_i},S_{i+1})$.

Let $i\geq0$ be fixed. We calculate first $P(T_0^{i+1}\leq t, \xi
^{i+1}=\xi|X^{S_i},S_{i+1})$ for $t\geq0, \xi\in\mathcal{P}$.
The theorem below shows that the location where $S_{i+1}$ first appears
is uniformly distributed on $\mathcal{T}_{S_i}$.
%
%th4 #&#
\begin{theorem}\label{distT0}
For any $t\geq0$, $\xi\in\mathcal{P}$, we have
\[
P\bigl(T_0^{i+1}\leq t, \xi^{i+1}=
\xi|X^{S_i},S_{i+1}\bigr)=\lambda\bigl(\bigl\{u\dvtx u\leq t, u<
\beta_{S_i}, \xi\in\mathcal{T}_{S_i}(u)\bigr\}
\bigr)/L_{S_i},
\]
where $\lambda$ is the Lebesegue measure and $\beta_{S_i}:=\inf\{t\dvtx |X^{S_i}|=1 \}$.
\end{theorem}
Theorem~\ref{distT0} is proved in Section~\ref{proofdistT0}.

For fixed $j\geq0$, with much more complicated argument and
discussion, we can calculate the conditional distribution
$P(T_{j+1}^{i+1}\in
B, \xi_{j+1}^{i+1}=\vec{\xi} | X^{S_i},S_{i+1},\break  T_0^{i+1},
\xi^{i+1},\ldots, T_j^{i+1},\xi_j^{i+1})$ for arbitrary
$B\in\mathcal{B}(R^{+})$ and $\vec{\xi}\in\mathcal{P}^{i+2}$. The
detailed discussion is given in the supplemental article \cite{supplement}.
The corresponding results are divided into 3 cases and stated below.

\textit{Case} 1: $\xi_j^{i+1}=\xi_0^{i+1}$. Since $\pi_{[0,i]}(\xi_0^{i+1})=\vec
{\varnothing}$, and the time point at
which $S_{i+1}$ appears is unique [cf. Proposition~\ref{Sgt}(iv)], in this
case the next event at time point $T_{j+1}^{i+1}$ must be a coalescence.
We have the following theorem.
%
%th5 #&#
\begin{theorem}\label{distTxi1} Suppose that $\xi_j^{i+1}=\xi_0^{i+1}$.

\textup{(i)} For $\vec{\xi}\in\mathcal{P}^{i+2}$ satisfying $\pi
_{i+1}(\vec{\xi})\neq\pi_{i+1}(\xi^{i+1}_{0})\cup
\pi_{i}(\vec{\xi})$, we
have
\[
P\bigl(\xi_{j+1}^{i+1}=\vec{\xi} | X^{S_i},S_{i+1},
T_0^{i+1},\xi _0^{i+1},
\ldots,T_j^{i+1},\xi_j^{i+1}\bigr)=0.
\]

\textup{(ii)} For $\vec{\xi}\in\mathcal{P}^{i+2}$ satisfying $\pi_{i+1}(\vec{\xi
})=\pi_{i+1}(\xi^{i+1}_{0})\cup
\pi_{i}(\vec{\xi})$, we have for arbitrary
$B\in\mathcal{B}(R^{+})$
\begin{eqnarray*}
&&P\bigl(T_{j+1}^{i+1}\in B, \xi_{j+1}^{i+1}=
\vec{\xi} | X^{S_i},S_{i+1}, T_0^{i+1},\xi
_0^{i+1},\ldots,T_j^{i+1},
\xi_j^{i+1}\bigr)\\
&&\qquad=\int_{T_j^{i+1}}^{\infty} I_{B}(t_{j+1})I_{\{t_{j+1}\dvtx \pi_{[0,i]}(\vec
{\xi})\in
\mathfrak{S}(X^{S_i}(t_{j+1}))\}}(t_{j+1})\\
&&\qquad\quad{}\cdot\exp\biggl\{-\int_{T_j^{i+1}}^{t_{j+1}}\bigl|X^{S_i}(t)\bigr|\,dt
\biggr\} \,dt_{j+1}.
\end{eqnarray*}
\end{theorem}
\begin{pf}
See Theorem B.10 in the supplemental article \cite{supplement}.
\end{pf}
\textit{Case} 2: $\xi_j^{i+1}\neq\xi_0^{i+1}$ and $\pi_i(\xi_j^{i+1})\neq
\varnothing$.
Because the time point at which $S_{i+1}$ appears is unique, in this case
$T_{j+1}^{i+1}$ must be a jump time of $X^{S_i}$. We define
%
%e3.17 #&#
\begin{equation}
\label{H} \mathcal{H}:=\inf\bigl\{t>T_j^{i+1}\dvtx
\pi_{[0,i]}\bigl(\xi_j^{i+1}\bigr)\notin
\mathfrak{S}\bigl(X^{S_i}(t)\bigr)\bigr\}.
\end{equation}
%
%th6 #&#
\begin{theorem}\label{distTxi3} Suppose that $\xi_j^{i+1}\neq\xi
_0^{i+1}$ and $\pi_i(\xi_j^{i+1})\neq\varnothing$. Then
for $\vec{\xi}\in\mathcal{P}^{i+2}$ satisfying $\pi_{i+1}(\vec{\xi
})=\pi_{i+1}(\xi^{i+1}_{0})\cup
\pi_{i}(\vec{\xi})$, $\pi_{i}(\xi_j^{i+1})\subset\pi_{i}(\vec{\xi})$
and $\pi_{[0,i]}(\vec{\xi})\in\mathfrak{S}(X^{S_i}(\mathcal{H}))$, it
holds that
\[
P\bigl(T_{j+1}^{i+1}=\mathcal{H}, \xi_{j+1}^{i+1}=
\vec{\xi } | X^{S_i},S_{i+1}, T_0^{i+1},
\xi_0^{i+1},\ldots,T_j^{i+1},
\xi_j^{i+1}\bigr)=1,
\]
where $\mathcal{H}$ is defined by (\ref{H}).
\end{theorem}
\begin{pf}
See Theorem B.11 in the supplemental article \cite{supplement}.
\end{pf}
\textit{Case} 3: $\xi_j^{i+1}\neq\xi_0^{i+1}$ and $\pi_i(\xi_j^{i+1})=\varnothing$.
In this case there is a potential recombination which generates again a
new lineage carrying $\xi_0^{i+1}$. Let $\mathcal{H}$ be defined by~(\ref{H}). If the waiting time is smaller than $\mathcal{H}-T_j^{i+1}$,
then recombination happens; otherwise no
recombination will happen and the lineage which carrys $\xi_j^{i+1}$
will follow the change of $X^{S_i}$.
In what follows, for an arbitrary $\vec{\xi}\in\mathcal{P}^{i+2}$,
we define
%
%e3.18 #&#
\begin{equation}
\label{hxi}h(\vec{\xi}):=\min\bigl\{l\dvtx \pi_p(\vec{\xi })=
\varnothing, \mbox{ for all } l< p\leq i\bigr\},
\end{equation}
if $\pi_i(\vec{\xi})=\varnothing$, otherwise
we set $h(\vec{\xi}):= i$.

%th7 #&#
\begin{theorem}\label{distTxi2} Suppose that $\xi_j^{i+1}\neq\xi
_0^{i+1}$ and $\pi_i(\xi_j^{i+1})=\varnothing$.

\textup{(i)} For $\vec{\xi}=\xi_0^{i+1}$ and
$B\in\mathcal{B}(R^{+})$, we
have
\begin{eqnarray*}
&&P\bigl(T_{j+1}^{i+1}\in B, \xi_{j+1}^{i+1}=
\xi_0^{i+1} | X^{S_i},S_{i+1},
T_0^{i+1},\xi _0^{i+1},
\ldots,T_j^{i+1},\xi_j^{i+1}\bigr)
\\
&&\qquad=\int_{T_j^{i+1}}^{\mathcal{H}} I_{B}(t_{j+1})
\biggl(\int_{S_{h(\xi
_j^{i+1})+1}}^{S_{i+1}}2^{-1}\rho p(v)\,dv
\biggr)
\\
&&\qquad\quad{}\cdot \exp \biggl\{-\bigl(t_{j+1}-T_j^{i+1}
\bigr)\int_{S_{h(\xi
_j^{i+1})+1}}^{S_{i+1}}2^{-1}\rho p(v)\,dv
\biggr\}\,dt_{j+1},
\end{eqnarray*}
where $h(\xi_j^{i+1})$ is specified by (\ref{hxi}).

\textup{(ii)} For $\vec{\xi}\in\mathcal{P}^{i+2}$ satisfying $\pi_{i+1}(\vec{\xi
})=\pi_{i+1}(\xi^{i+1}_{0})\cup
\pi_{i}(\vec{\xi})$, $\pi_{i}(\xi_j^{i+1})\subset\pi_{i}(\vec{\xi})$
and $\pi_{[0,i]}(\vec{\xi})\in\mathfrak{S}(X^{S_i}(\mathcal{H}))$, we have
\begin{eqnarray*}
&&P\bigl(T_{j+1}^{i+1}=\mathcal{H}, \xi_{j+1}^{i+1}=
\vec{\xi } | X^{S_i},S_{i+1}, T_0^{i+1},
\xi_0^{i+1},\ldots,T_j^{i+1},\xi
_j^{i+1}\bigr)
\\
&&\qquad=\exp \biggl\{-\bigl(\mathcal{H}-T_j^{i+1}\bigr)\int
_{S_{h(\xi
_j^{i+1})+1}}^{S_{i+1}}2^{-1}\rho p(v)\,dv \biggr\}.
\end{eqnarray*}

\textup{(iii)} For $\vec{\xi}\in\mathcal{P}^{i+2}$ satisfying neither \textup{(i)} nor
\textup{(ii)}, we have
\[
P\bigl(\xi_{j+1}^{i+1}=\vec{\xi} | X^{S_i},S_{i+1},
T_0^{i+1},\xi _0^{i+1},
\ldots,T_j^{i+1},\xi_j^{i+1}\bigr)=0.
\]
\end{theorem}
\begin{pf}
See Theorem B.12 in the supplemental article \cite{supplement}.
\end{pf}

%s3.2 #&#
\subsection{Spatial algorithm}
\label{algorithm}
Based on the random sequence $\{(S_i, Z^i), i\ge0\}$ discussed above,
in this subsection we describe our model of spatially moving algorithm
and study its statistical property.

%s3.2.1 #&#
\subsubsection{SC algorithm}

In this subsection, we describe briefly our spatially moving algorithm;
for details, see \cite{wang}. Our new algorithm is called $\mathit{SC}$, which
will recursively construct part graph $X^{S_i}$ with each branch
assigned some label $k\leq i$. All the branches with label $i$
form the local tree $\mathcal{T}_{S_i}$.

\textit{Step} 1. Construct a standard Kingman's coalescent tree (cf. \cite
{Kingman}) $\mathcal{T}_0$
at the position $S_0=0$ (the left end point
of the sequence) and assign each branch of the tree with the label $0$.
Let $X^{0}=\mathcal{T}_0$.

\textit{Step} 2. Assume that we have already constructed $X^{S_i}$ along with
local tree $\mathcal{T}_{S_i}$. Take the next recombination point
$S_{i+1}$ along the sequence according to the distribution
\[
P\bigl(S_{i+1}> s|X^{S_i}\bigr)=\exp\biggl\{-\rho
L_{S_i}\bigl(X^{S_i}\bigr)\int_{S_i}^{s\vee S_i}2^{-1}p(r)\,dr
\biggr\}.
\]
If
$S_{i+1}\geq1$, stop; otherwise, go to step 3.\vadjust{\goodbreak}

\textit{Step} 3. Uniformly choose a recombination location on
$\mathcal{T}_{S_i}$. For $j=0$, let $T_j^{i+1}$ denote the latitude
(i.e., the height from the bottom to the location) of the chosen
location.

\textit{Step} 4. At the recombination location, a new branch with label $i+1$
is created by forking off the recombination node and moving backward
in time (i.e., along the direction of increasing latitude). With
equal exponential rate 1, the new branch will have a tendency to
coalesce to each branch in $X^{S_i}$ which has higher latitude than
$T_j^{i+1}$. Thus, if there are $l$ branches in $X^{S_i}$ at the
current latitude, then the waiting time before coalescence is
exponentially distributed with parameter $l$. Note at different
latitude there may be different number $l$ of branches. Let the branch
to which the new branch coalesces be called EDGE, and let
$T_{j+1}^{i+1}$ be the latitude of the coalescent point and regard
$j+1$ as $j$ in the next step.

\textit{Step} 5. If the EDGE is labeled with $i$, go to step 6; if the EDGE
is labeled with some $k$ less than $i$, then a potential
recombination event should be considered. The waiting time $t$ of
the possible recombination event on the EDGE is exponentially
distributed with parameter $\int_{S_{k+1}}^{S_{i+1}}2^{-1}\rho
p(u)\,du$.

\begin{itemize}
\item \textit{Case} 5.1. If $T_{j}^{i+1}+t$ is less than the latitude of the
upper node of the EDGE
which is denoted by $\mathcal{H}$, then it is the next recombination location.
Let $T_{j+1}^{i+1}=T_{j}^{i+1}+t$, the part of the branch above
$T_{j+1}^{i+1}$ is no longer called EDGE. Regard $j+1$ as $j$ and go to
step 4.

\item \textit{Case} 5.2. If $T_{j}^{i+1}+t\geq\mathcal{H}$, choose the upper
edge of the current EDGE with the larger label to
be the next EDGE.
Let $T_{j+1}^{i+1}=\mathcal{H}$, regard $j+1$ as $j$ and go to step 5.
\end{itemize}

\textit{Step} 6. Let $X^{S_{i+1}}$ be the collection of all the branches in
$X^{S_{i}}$ and all the new branches labeled $i+1$. Starting from
each node $1\leq m\leq N$ at the bottom of the graph, specify a
path moving along the edges in $X^{S_{i+1}}$ increasing latitude,
until the top of the graph. Whenever a recombination node is
encountered, choose the edge with the larger label. The collection of
all the paths then form the local tree $\mathcal{T}_{S_{i+1}}$.
Update all the branches in $\mathcal{T}_{S_{i+1}}$ with label $i+1$.
%s3.2.2 #&#
\subsubsection{Some explanation of SC}
Assume that we have already constructed $X^{S_i}$.
Then the local tree $\mathcal{T}_{S_i}$ and the breakpoints
$S_1,\ldots,S_i$ are all known. Thus, steps~2 and~3 are feasible. Moreover, once we have constructed
$X^{S_i}$, then the ancestral material of each edge
in $X^{S_i}$, expressed as an $(i+1)$-dimensional
$\mathcal{P}$-valued vector [cf. (\ref{f^Si})], is also implicitly
known. In step 4, if we denote by
$\vec{z}=(z_0,z_1,\ldots,z_i)\in\mathcal{P}^{i+1}$ the ancestral
material of the edge where the recombination location lies, then it is
implicitly assumed that the ancestral material carried on the new
branch with label $i+1$ is the $(i+2)$-dimensional $\mathcal{P}$-valued vector
$(\varnothing,\ldots,\varnothing,z_i)$. If for $j=0$ we write $\xi^{i+1}$
for $z_i$, and denote
$\xi_0^{i+1}:=(\varnothing,\ldots,\varnothing,\xi^{i+1})$, then in
steps 4 and~5 the algorithm specifies a path describing how the
ancestral material $\xi_0^{i+1}$ coalesces to
$X^{S_i}$ by coalescence and leaves
$X^{S_i}$ by recombination.
When the EDGE in step 4
is labeled $i$, then it is implicitly assumed that the path
carrying the ancestral material $\xi_0^{i+1}$ extends continuously
along the edges with the larger label in $X^{S_{i}}$, starting from the
EDGE until the top of the graph. In step 6, the algorithm formulates
$X^{S_{i+1}}$ with the branches in $X^{S_i}$ and the new branches
created in step 4. For an edge in $X^{S_i}$ carrying an ancestral material
$(z_0,z_1,\ldots,z_i)\in\mathcal{P}^{i+1}$, when the edge is viewed as
an edge in $X^{S_{i+1}}$, it is implicitly assumed that its ancestral
material is updated to $(z_0,z_1,\ldots,z_i,z_{i+1})\in\mathcal
{P}^{i+2}$ by the following rule: (i) if the edge is on the path of
$\xi_0^{i+1}$ specified above, then $z_{i+1}=z_i \cup\xi^{i+1}$; (ii)
if the edge is not on the path of $\xi_0^{i+1}$, then $z_{i+1}=z_i$ on
the part of the edge below the latitude $T_0^{i+1}$, and $z_{i+1}=z_i
\setminus\xi^{i+1}$ on the part of the edge above the latitude $T_0^{i+1}$.
%s3.2.3 #&#
\subsubsection{Distribution of the ARG generated by SC}
In this subsection, we shall show that the probability distribution
of the ARG generated by the $\mathit{SC}$ algorithm coincides with that generated
by the back-in-time algorithm as specified in Theorem~\ref{space}.
To this end, we denote by $|Bp|$ the maximum $i$ such that
$S_i < 1$. For each $0\leq i < |Bp|$, we denote by $Z^{i+1}(t)$ the
ancestral material [represented as an $(i+2)$-dimensional $\mathcal
{P}$-valued vector] at the latitude $t\geq T_0^{i+1}$ on the
path of edges carrying $\xi_0^{i+1}$, and set $Z^{i+1}(t) = \vec
{\varnothing}$ for $t < T_0^{i+1}$, with $\vec{\varnothing}$ representing
the null vector in $\mathcal{P}^{i+2}$. For $i\geq|Bp|$, we set
$S_{i+1}=1,  T_0^{i+1}=\infty$, and $Z^{i+1}(t)\equiv\vec{\varnothing}$.
Further, we write $Z^0(t)=\mathcal{T}_0(t)$ for $t\geq0$.
%
%pr16 #&#
\begin{proposition}\label{finitedim}
With the above convention, the finite dimensional distribution of the
random sequence $\{(S_i, Z^i),i\geq0\}$ generated by the $\mathit{SC}$
algorithm is the same as that developed in Section~\ref{sequence}.
\end{proposition}
Proposition~\ref{finitedim} is proved in Section~\ref{prooffinitedim}.

By virtue of the above proposition, we are in a position to prove the
following most important theorem of this paper.
%
%th8 #&#
\begin{theorem}\label{coincide}
Let $(G,\mathcal{B}(G),P)$ be the probability space specified in
Theorem~\ref{space}, and denote by $\tilde{P}$ the probability
distribution on $G$ generated by $\mathit{SC}$ algorithm. Then we have $\tilde{P}=P$.
\end{theorem}
Theorem~\ref{coincide} is proved in Section~\ref{proofcoincide}.

%s4 #&#
\section{Proofs and technical lemmas}\label{proofs}

%s4.1 #&#
\subsection{Proof of Proposition \texorpdfstring{\protect\ref{fil}}{14}}\label{prooffil}
The proof of Proposition~\ref{fil} needs the following lemma.\vadjust{\goodbreak}
%
%le1 #&#
\begin{lemma}\label{F^Si}
For each $i\geq1$, we have
%
%e4.1 #&#
\begin{equation}
\label{Fsi} \mathcal{F}^{S_i}=\sigma\bigl(X^0,S_1,
\ldots,S_i;X^{S_1},\ldots,X^{S_i}\bigr) =\sigma
\bigl(X^{S_i}\bigr).
\end{equation}
\end{lemma}
\begin{pf}
By the classical theory of jump processes (cf., e.g., \cite{Yan92},
Definition~11.48 and Corollary~5.57), we have
%
%e4.2 #&#
\begin{equation}
\label{F^s}\mathcal{F}^s=\bigcup_{i=0}^{\infty}
\bigl(\mathcal {G}_i\cap\bigl\{S_i^{\prime}\leq
s<S_{i+1}^{\prime}\bigr\} \bigr)
\end{equation}
and
\[
\mathcal{F}^{S_i^{\prime}}=\mathcal{G}_i,
\]
where $\mathcal{G}_i=\sigma(X^0,S_1^{\prime},\ldots,S_i^{\prime};X^{S_1^{\prime}},\ldots,X^{S_i^{\prime}})$.
Since $X^{s}=X^{1}$ for $s\geq1$, one can check that $\mathcal
{F}^{1}=\mathcal{F}^{\infty}$, $S_i=S_i^{\prime}\wedge1$, $\sigma
(S_i)=\sigma(S_i^{\prime})$ and $X^{S_i}=X^{S_i^{\prime}}$ for $i\geq1$. Hence,
we have $\mathcal{G}_i=\sigma(X^0,S_1,\ldots,S_i;X^{S_1},\ldots,X^{S_i})$.
Therefore, we have
\[
\mathcal{F}^{S_i}=\mathcal{F}^{S_i^{\prime}}\cap\mathcal
{F}^{1}=\mathcal{F}^{S_i^{\prime}} =\sigma\bigl(X^0,S_1,
\ldots,S_i;X^{S_1},\ldots,X^{S_i}\bigr).
\]

Next, for an arbitrary $g\in G$, suppose that $g(t)$ is expressed as
$g(t)= \{f_1,f_2,\ldots,f_{k}\}$,
then by Proposition~\ref{Sgt}(v) we can show that
%
%e4.3 #&#
\begin{equation}
\label{f^Si} f^{S_i}_j=\sum_{l=0}^{i-1}
\pi_l(\vec{z}_j)I_{[S_l,S_{l+1})}+ \pi_i(
\vec {z}_j)I_{[S_i,1)}.
\end{equation}
Consequently, for $l\leq i$, we have $X^{S_l}=\pi
_{[0,S_l]}^{G}(X^{S_i})$ and $S_l (g)= S_l(X^{S_i}(g))$,
and the second equality of (\ref{Fsi}) follows.
\end{pf}

\begin{pf*}{Proof of Proposition~\ref{fil}}
By (\ref{f^Si}), we have
$Z^l(g)=Z^l(X^{S_i}(g))$ for all $l\leq i$, hence
\[
\sigma\bigl(S_1,\ldots,S_i;Z^0,Z^1,
\ldots,Z^i\bigr)\subset \sigma\bigl(X^0,S_1,
\ldots,S_i;X^{S_1},\ldots,X^{S_i}\bigr)=\sigma
\bigl(X^{S_i}\bigr).
\]
To show the inverse inclusion, we put
%
%e4.4 #&#
\begin{equation}
\label{Omigai} \Omega_i:=[0,1]^i \times
\mathcal{S}_{[0,\infty)}(\mathcal{R}) \times \prod
_{l=1}^{i} \mathcal{S}_{[0,\infty)}\bigl(
\mathcal{P}^{l+1}\bigr),
\end{equation}
where $\mathcal{R}$ is the collection of all the partitions of $\{
1,2,\ldots,N\}$, and $S_{[0,\infty)}(\mathcal{P}^{l+1})$ [resp., $\mathcal
{S}_{[0,\infty)}(\mathcal{R})$] equipped with the Skorohod\vadjust{\goodbreak} topology are
the spaces of all the
$\mathcal{P}^{l+1}$-valued (resp.,\vadjust{\goodbreak} $\mathcal{R}$-valued) right
continuous piecewise constant functions
on $[0,\infty)$ with at most finitely many discontinuity points.
Define
%
%e4.5 #&#
\begin{equation}
\label{Phii} \Phi_i:=\bigl(S_1,\ldots,S_i;Z^0,Z^1,
\ldots,Z^i\bigr).
\end{equation}
From
Propositions \ref{3Gneighborhood} and~\ref{Zjcon}, we see
that $\Phi_i$ is a continuous map from $G$ to the Polish space $\Omega
_i$. Denote by $\mathcal{H}_i=\Phi_i(G)$. By (\ref{f^Si}), we have
$\mathcal{H}_i=\Phi_i(X^{S_i}(G))$. By (\ref{f^Si}) and Proposition~\ref
{V} one can check that $\Phi_i$ restricted to $X^{S_i}(G)$ is an
injective map. Below we write $G_i:=X^{S_i}(G)$. Note that $G_i=\{g\in
G\dvtx g=X^{S_i}(g)\}$ is a Borel subset of the Polish space $G$. Hence
$\mathcal{H}_i=\Phi_i(G_i)$ is a Borel subset of $\Omega_i$ and
$(\Phi_i|_{G_i})^{-1}\dvtx \mathcal{H}_i \mapsto G_i$ is Borel measurable
(cf. \cite{Cohn80}, Theorems~8.3.5 and 8.3.7). Define a map $\Upsilon_i\dvtx \Omega_i\mapsto G$ by setting $\Upsilon_i(\omega)=(\Phi
_i|_{G_i})^{-1}(\omega)$ if $\omega\in\mathcal{H}_i$ and $\Upsilon
_i(\omega)=g_0$ if $\omega\notin\mathcal{H}_i$, where $g_0$ is a
fixed element in $G$. Since $(\Phi_i|_{G_i})^{-1}$ is Borel measurable
and $\mathcal{H}_i$ is a Borel subset of $\Omega_i$, hence $\Upsilon_i$
is also Borel measurable.
Noticing that $X^{S_i}(g)=\Upsilon_i(\Phi_i(g))$, we conclude that
$\sigma(X^{S_i})\subset\sigma(S_1,\ldots,S_i;Z^0,Z^1,\ldots,Z^i)$,
completing the proof.
\end{pf*}

%s4.2 #&#
\subsection{Proof of Theorem \texorpdfstring{\protect\ref{distSi}}{2}}\label{proofdistSi}
The proof of Theorem~\ref{distSi} requires the following two lemmas.
%
%le2 #&#
\begin{lemma}\label{Si1}
For $0\leq u<s<1$, it holds that
\[
P\bigl(\bigl\{g\dvtx \pi_{[0,s]}^G(g)=\pi_{[0,u]}^G(g)
\bigr\}|\mathcal{F}^u\bigr)=\exp\biggl\{-\rho L_u
\bigl(X^u\bigr)\int_{u}^{s}2^{-1}p(r)\,dr
\biggr\}.
\]
\end{lemma}
\begin{pf}
For $s\geq0$, we set $\gamma^s:=\gamma^s(g)=\inf\{n\dvtx \tau
_{n+1}^s(g)=\infty\}$. Note that $X^s(\cdot)(g):=\pi_{[0,s]}^G(g)$ is
uniquely determined by $\{X^s(0),\tau_1^{s}, X^s(\tau_1^{s}),\tau
_2^{s},\break X^s(\tau_2^{s}),
\ldots\}$. Therefore
for $0\leq u<s$, $\pi_{[0,s]}^G(g)=\pi_{[0,u]}^G(g)$ if and only if
$\tau_n^{s}=\tau_n^{u}$ and $X^s(\tau_n^{s})=X^u(\tau_n^{u})$ for all
$1\leq n \leq\gamma^u$, and $\gamma^s=\gamma^u$. Write $G^u=X^u(G)$.
Since $G^u=\{g\in G\dvtx  g=\pi_{[0,u]}^G(g)\}$, hence $G^u\in\mathcal
{B}(G)$. For $k\geq1$, we set $A_k=\{g\in G^u\dvtx   \gamma^u(g^u)=k\}$
and define
\[
F_k(g):=\bigl(\tau_1^{u}(g),X^u
\bigl(\tau_1^{u}\bigr) (g),\ldots,\tau_k^{u}(g),
X^u\bigl(\tau_k^{u}\bigr) (g)\bigr)
\]
for $g\in A_k$. Then $F_k$ is a measurable map from $A_k$ to
$B_k:=F_k(A_k)\subset(R^{+}\times E^u)^k$. Since $A_k$ is a Borel
subset of the Polish space $G$ and $F_k$ is injective, $B_k$ is a Borel
subset of $(R^{+}\times E^u)^k$ (cf., e.g., \cite{Cohn80}, Theorem~8.3.7). By the one to one correspondence between $A_k$ and $B_k$, one
can check that if $g\in G$ satisfies $(\tau_1^{s}(g),X^s(\tau
_1^{s})(g),\ldots,\tau_k^{s}(g),
X^s(\tau_k^{s})(g)) \in B_k$, then $\pi_{[0,s]}^G(g)=\pi_{[0,u]}^G(g)$.

For $A\in\mathcal{B}(G^u), A\subset A_k$, denote by $B=F_k(A)\subset
B_k$. Applying Proposition~\ref{distribution2}, we have
\begin{eqnarray*}
&&P\bigl(\bigl\{g\dvtx X^u\in A,\pi_{[0,s]}^G(g)=
\pi_{[0,u]}^G(g)\bigr\}\bigr)
\\
&&\qquad=P\bigl(\bigl\{g\dvtx \bigl(X^u(0),\tau_1^{u},X^u
\bigl(\tau_1^{u}\bigr),\ldots,\tau_k^{u},X^u
\bigl(\tau _k^{u}\bigr)\bigr)\in B,\pi_{[0,s]}^G(g)=
\pi_{[0,u]}^G(g)\bigr\}\bigr)
\\
&&\qquad=P\bigl(\bigl\{g\dvtx \bigl(X^s(0),\tau_1^{s},X^s
\bigl(\tau_1^{s}\bigr),\ldots,\tau_k^{s},X^s
\bigl(\tau _k^{s}\bigr)\bigr)\in B\bigr\}\bigr)
\\
&&\qquad=\int_{R^+\times E}\cdots\int_{R^+\times E}\cdots \int
_{R^+\times E}I_{B}(y_0,t_1,y_1,
\ldots,t_k,y_k)
\\
&&\qquad\quad{} \cdot\exp\Biggl\{-\sum_{n=0}^{k-1}q
\bigl(y_{n},\bigl[y_{n}^{s}\bigr]^{c}
\bigr) (t_{n+1}-t_{n})\Biggr\}
\\
&&\qquad\quad dt_{k+1}q(y_{k},dy_{k+1})\cdots
dt_{n+1}q(y_{n},dy_{n+1})\cdots
dt_1q(y_0,dy_1).
\end{eqnarray*}
For $\omega=(t_1,y_1,\ldots,t_k,y_k)\in B_k$, if we set $g_\omega:=y_0I_{[0,t_1)}(t)+\sum_{j=1}^{k}y_jI_{[t_j,t_{j+1})}(t)$, with the
convention that $t_{k+1}=\infty$, then $g_\omega\in G^u$ and it is easy
to check that $L_u(g_\omega)=\sum_{n=0}^{k-1}|\pi
_u(y_{n})|(t_{n+1}-t_{n})$. Since $y_n\in E^{u}$, by (\ref{qxA}) we can
derive that for $u< s$, $q(y_n,[y_n^{s}]^{c})=q(y_n,[y_n^{u}]^{c})+|\pi
_u(y_n)|\int_u^s 2^{-1}\rho p(r)\,dr$. Therefore, from the above we get
\begin{eqnarray*}
&&P\bigl(\bigl\{g\dvtx X^u\in A,\pi_{[0,s]}^G(g)=
\pi_{[0,u]}^G(g)\bigr\}\bigr)
\\
&&\qquad=\int_{R^+\times E}\cdots\int_{R^+\times E}\cdots \int
_{R^+\times E}I_{B}(y_0,t_1,y_1,
\ldots,t_k,y_k)
\\
&&\qquad\quad{} \cdot\exp\Biggl\{-\sum_{n=0}^{k-1}q
\bigl(y_{n},\bigl[y_{n}^{u}\bigr]^{c}
\bigr) (t_{n+1}-t_{n})\Biggr\}
\\
&&\qquad\quad{} \cdot\exp\Biggl\{-\sum_{n=0}^{k-1}\bigl|
\pi_u(y_n)\bigr|(t_{n+1}-t_n) \int
_{u}^{s}2^{-1}\rho p(r)\,dr\Biggr\}
\\
&& \qquad\quad dt_{k+1}q(y_{k},dy_{k+1})\cdots
dt_{n+1}q(y_{n},dy_{n+1})\cdots
dt_1q(y_0,dy_1)
\\
&&\qquad= \int_G I_{\{X^u\in A\}}\exp\biggl\{-\rho
L_u\bigl(X^u\bigr)\int_{u}^{s}2^{-1}p(r)\,dr
\biggr\}P(dg).
\end{eqnarray*}
The proof is completed by noticing that $G=\sum_k A_k$ and
$\mathcal{F}^u=\sigma(X^u)$.
\end{pf}

%le3 #&#
\begin{lemma}\label{Si2}
For an arbitrary $(\mathcal{F}^s)$-stopping time $\tau$, we write
$T_{\tau}=\break \inf\{s>\tau: X^s\neq X^{\tau}\}$, then the distribution
of $T_{\tau}$ conditioning on $\mathcal{F}^{\tau}$ is: for $s<1$,
\[
P\bigl(T_{\tau}> s|\mathcal{F}^{\tau}\bigr)= \exp\biggl\{ -\rho
L_{\tau}\bigl(X^{\tau}\bigr)\int_{\tau}^{s\vee\tau}2^{-1}p(r)\,dr
\biggr\}.
\]
\end{lemma}

\begin{pf}
For $n\geq1$, we define
\[
\tau^{(n)}= \sum_{j=1}^\infty
\frac{j}{2^n}I_{\{{(j-1)}/{2^n}\leq
\tau<{j}/{2^n}\}}+\infty I_{\{\tau=\infty\}}.
\]
Then $\{\tau^{(n)}\}$ are countably valued stopping times
and $\tau^{(n)}\downarrow\tau$. By the convergence theorem of
conditional expectations
(cf., e.g., \cite{Yan92}, Theorem~2.21),\vadjust{\goodbreak} for $s<1$ employing Lemma~\ref
{Si1} we obtain
\begin{eqnarray*}
P\bigl(T_{\tau}> s|\mathcal{F}^{\tau}\bigr)I_{\{\tau< s\}}& =&P
\bigl(T_{\tau}> s,\tau< s|\mathcal{F}^{\tau}\bigr)
\\
& =&\lim_{n\rightarrow\infty}P\bigl(T_{\tau}> s,\tau^{(n)} <
s|\mathcal{F}^{\tau
^{(n)}}\bigr)
\\
&=&\lim_{n\rightarrow\infty}\sum_{j=1}^{\infty}P
\bigl(T_{\tau}> s,j/2^n < s|\mathcal{F}^{j/2^n}
\bigr)I_{\{\tau^{(n)}=j/2^n\}}
\\
&=&\lim_{n\rightarrow\infty}\sum_{j/2^n < s}\exp
\biggl(-\rho L_{j/2^n}\int_{j/2^n}^{s}2^{-1}p(r)\,dr
\biggr)I_{\{\tau^{(n)}=j/2^n\}}
\\
&=&\lim_{n\rightarrow\infty}\exp\biggl(-\rho L_{\tau^{(n)}}\int
_{\tau^{(n)}}^{s}2^{-1}p(r)\,dr
\biggr)I_{\{\tau^{(n)}< s \}}
\\
&=&\exp\biggl(-\rho L_{\tau}\int_{\tau}^{s}2^{-1}p(r)\,dr
\biggr)I_{\{\tau< s\}}.
\end{eqnarray*}
On the other hand, we have
\begin{eqnarray*}
P\bigl(T_{\tau}> s|\mathcal{F}^{\tau}\bigr)I_{\{\tau\geq s\}} &=&P
\bigl(T_{\tau}> s,\tau\geq s|\mathcal{F}^{\tau}\bigr)=P\bigl(\tau
\geq s|\mathcal{F}^{\tau}\bigr)
\\
&=&I_{\{\tau\geq s\}}=\exp\biggl\{-\rho L_{\tau}\int_{\tau}^{s\vee\tau}2^{-1}p(r)\,dr
\biggr\}I_{\{\tau\geq
s\}}.
\end{eqnarray*}
\upqed\end{pf}

\begin{pf*}{Proof of Theorem~\ref{distSi}}
If $\tau=S_i$, then we have $S_{i+1}=T_\tau\wedge1$. Therefore, the
theorem is a direct consequence of the above lemma.
\end{pf*}

%s4.3 #&#
\subsection{Proof of Theorem \texorpdfstring{\protect\ref{distZ^0}}{3}}\label{proofdistZ^0}

Comparing (\ref{Z^O}) and (\ref{Ts}), it is apparent that
$Z^0(t)=\mathcal{T}_0(t)$. In Proposition~\ref{distribution2}, let
$s=u=0$ and $k=0$, then for $B\in\mathcal{B}( (R^{+}\times E^0)^{n})$
we have
\begin{eqnarray*}
&&P\bigl\{\bigl(\tau_1^{0},X^0\bigl(
\tau_1^{0}\bigr),\ldots,\tau _n^{0},X^0
\bigl(\tau_n^{0}\bigr)\bigr)\in B\bigr\}\\
&&\qquad=\int_0^{\infty}\,dt_1\cdots\int
_0^{\infty}\,dt_{n}\int_Eq(y_0,dy_1)
\cdots \int_E q(y_{n-1},dy_n)
\\
&&\qquad\quad{} \cdot I_{\{t_1<\cdots< t_{n}\}}(t_1,\ldots,t_{n})\prod
_{j=0}^{n-1}I_{E^{0}}(y_{j+1})\\
&&\qquad\quad{}\cdot\exp
\Biggl\{-\sum_{j=0}^{n-1}q
\bigl(y_j,\bigl[y_j^0\bigr]^{c}
\bigr) (t_{j+1}-t_j)\Biggr\}.
\end{eqnarray*}
Note that $q(y,[y^0]^{c})=|\pi_0^E(y)|(|\pi_0^E(y)|-1)/2$ for $y\in
E^0$. Identifying $y$ with $\pi_0^E(y)$ for $y\in E^0$, we see that $\{
\mathcal{T}_0(t), t\geq0\}$ is a standard Kingman's coalescent tree.

%s4.4 #&#
\subsection{Proof of Theorem \texorpdfstring{\protect\ref{distT0}}{4}}\label{proofdistT0}
We first make some preparing discussions.
For $ 0\leq j \leq k, a\in[0,1)$ and $\xi\in\mathcal{P}$, we put
$A_{k,a,j,\xi}=\{ \gamma^a=k,S_i\leq a, \xi\in\mathcal{T}_{a}(\tau
_{j}^{a})\}$. For notational convenience, below we use also
$A_{k,a,j,\xi}$ to denote the indicator function of $A_{k,a,j,\xi}$.
For $j\geq0$, $a\in[0,1)$, $s\in[0,1)$, $\varepsilon>0$ and $\xi\in
\mathcal{P}$, we define
%
%e4.6 #&#
\begin{eqnarray}
&&\label{K1}K_{a,s,\varepsilon,j,\xi}=\bigl\{s-\varepsilon\leq S_{i+1}< s,
\tau_{j+1}^{s}\leq t, \tau_{j}^a<
\tau_{j+1}^{s}< \tau _{j+1}^{a},
\nonumber
\\
&&\hspace*{58pt}\exists r\in (s-\varepsilon,s) \mbox{ and } f\in X^{s}\bigl(\tau
_{j+1}^{s}\bigr),\\
&&\hspace*{100pt} \mbox{s.t. } f(r)=\xi, f(u)=\varnothing \mbox{ for
all } u<r\bigr\}.\nonumber
\end{eqnarray}

%le4 #&#
\begin{lemma}\label{distZ1} Let $s\in[0,1)$ and $\varepsilon>0$ be such
that $a< s-\varepsilon$.
Then for any $t>0$, we have
\begin{eqnarray*}
&&P\bigl(K_{a,s,\varepsilon,j,\xi} |X^a\bigr)A_{k,a,j,\xi}
\\
&&\qquad=C(a,j,s,\varepsilon)P\bigl(\{s-\varepsilon\leq S_{i+1}< s\}
|X^a\bigr)A_{k,a,j,\xi},
\end{eqnarray*}
where
%
%e4.7 #&#
\begin{eqnarray}
\label{Cajs} &&C(a,j,s, \varepsilon)\nonumber\\
&&\qquad:= \biggl(1-\exp\biggl\{-\rho
L_a\int_{s-\varepsilon
}^{s}2^{-1}p(r)\,dr
\biggr\} \biggr)^{-1}\cdot\int_{s-\varepsilon}^{s} 2^{-1}\rho p(r)\,dr
\nonumber
\\[-8pt]
\\[-8pt]
\nonumber
&&\qquad\quad{}\cdot\int_{\tau_j^a\wedge t}^{\tau_{j+1}^a\wedge t}\exp
\Biggl\{- \Biggl[\bigl|
\mathcal{T}_{a}\bigl(\tau_j^a\bigr)\bigr|
\bigl(t^{\prime}-\tau_j^a\bigr)
+ \sum_{l=0}^{j-1}\bigl|\mathcal{T}_{a}
\bigl(\tau _l^a\bigr)\bigr|\bigl(\tau_{l+1}^{a}-
\tau_l^a\bigr) \Biggr] \\
&&\hspace*{195pt}\qquad\quad{}\cdot\int_{s-\varepsilon}^{s}2^{-1}
\rho p(r)\,dr \Biggr\}\,dt^{\prime}.
\nonumber
\end{eqnarray}
\end{lemma}
\begin{pf}
We need only to check the lemma in the case that $A_{k,a,j,\xi}\neq
\varnothing$. Note that $A_{k,a,j,\xi}\in\mathcal{F}^a=\sigma(X^a)$.
Take an arbitrary set $H\in\sigma(X^a)$. Define $B=F_k(H\cap
A_{k,a,j,\xi})$, where $F_k$ is specified as in the proof of Lemma~\ref{Si1}.
Then $B\subset(R^{+}\times E^a)^k$ and
\[
H\cap A_{k,a,j,\xi}=\bigl\{g\dvtx \bigl(\tau_1^{a}(g),X^a
\bigl(\tau_1^{a}\bigr) (g),\ldots,\tau_k^{a}(g),
X^a\bigl(\tau_k^{a}\bigr) (g)\bigr) \in B
\bigr\}.
\]
Suppose that $B\neq\varnothing$.
One can check that
%
%e4.8 #&#
\begin{eqnarray}
\label{inB} &&K_{a,s,\varepsilon,j,\xi}\cap H\cap A_{k,a,j,\xi} \nonumber\\
&&\qquad= \bigl\{g\dvtx
\bigl(\tau _1^{s},X^{s}\bigl(
\tau_1^{s}\bigr), \ldots,\tau_j^{s},X^{s}
\bigl(\tau_j^{s}\bigr),
\nonumber\\
&&\hspace*{22pt}\qquad\quad \tau _{j+1}^{s-\varepsilon}, X^{s-\varepsilon}\bigl(
\tau_{j+1}^{s-\varepsilon}\bigr),\ldots, \tau_{k}^{s-\varepsilon},
X^{s-\varepsilon}\bigl(\tau_{k}^{s-\varepsilon}\bigr)\bigr)\in B,
\\
&&\hspace*{5pt}\qquad\quad\tau_{j+1}^{s}\leq t, \tau_{j+1}^s<
\tau_{j+1}^{s-\varepsilon}, \exists r\in(s-\varepsilon,s) \mbox{ and } f\in
X^{s}\bigl(\tau_{j+1}^{s}\bigr),
\nonumber
\\
&&\hspace*{95pt}\qquad\quad \mbox{s.t. } f(r)=\xi, f(u)=\varnothing \mbox{ for all } u<r \bigr\}.
\nonumber
\end{eqnarray}
Because $B\subset F_k(\{\gamma^a=k, \xi\in\mathcal{T}_{a}(\tau
_{j}^{a})\})$, therefore, for
$(t_1,y_1,\ldots,t_j,y_j,\ldots,\break  t_k,y_k)\in B$, there exists $f_l \in
y_j$ such that $f_l(u)=\xi$ for $u\in[a,1)$. We set
\[
J_{y_j}=\bigl\{y\in E^s\dvtx y=R_{lu}(y_j),
u\in(s-\varepsilon,s) \mbox{ and } l \mbox { satisfies } f_l\in
y_j, f_l(u)=\xi\bigr\}
\]
and define
\begin{eqnarray*}
B^{\prime}&:=& \bigl\{\bigl(t_1,y_1,
\ldots,t_j,y_j,t^{\prime},y^{\prime},t_{j+1},
y_{j+1}, \ldots,t_k,y_k\bigr)\in
\bigl(R^{+}\times E^s\bigr)^{k+1}\dvtx
\\
&&\hspace*{6pt}(t_1,y_1,\ldots,t_j,y_j,
\ldots,t_k,y_k)\in B, t^{\prime}
\in(t_j,t_{j+1})\cap(0,t], y^{\prime}\in
J_{y_j} \bigr\}.
\end{eqnarray*}
With $\tau_{j+1}^{s},X^{s}(\tau_{j+1}^{s})$ in the place of
$t^{\prime},y^{\prime}$, we may write (\ref{inB}) as
\begin{eqnarray*}
&& K_{a,s,\varepsilon,j,\xi}\cap H\cap A_{k,a,j,\xi} \\
&&\qquad= \bigl\{g\dvtx \bigl(
\tau_1^{s},X^{s}\bigl(\tau_1^{s}
\bigr),\ldots,\tau_j^{s},X^{s}\bigl(\tau
_j^{s}\bigr),\tau_{j+1}^{s},
\\
&&\hspace*{33pt}\qquad X^{s}\bigl(\tau_{j+1}^{s}\bigr),
\tau_{j+1}^{s-\varepsilon}, X^{s-\varepsilon}\bigl(\tau_{j+1}^{s-\varepsilon}
\bigr),\ldots,\tau_{k}^{s-\varepsilon}, X^{s-\varepsilon}\bigl(
\tau_{k}^{s-\varepsilon}\bigr) \bigr) \in B^{\prime} \bigr\}.
\end{eqnarray*}
For $j+1\leq l\leq k$, we set $\vartheta_{l}:= \tau_{j+1}^{s}+\tau
_{l-j}^{s-\varepsilon}\circ\theta_{\tau_{j+1}^{s}}$. One can check that
$\vartheta_{l}=\tau_l^{s-\varepsilon}$ for each $l$. Employing the strong
Markov property and Proposition~\ref{distribution2}, we get
\begin{eqnarray*}
&&P(K_{a,s,\varepsilon,j,\xi}\cap H\cap A_{k,a,j,\xi
})\\
&&\qquad=\int_0^{\infty}\,dt_1\cdots\int
_0^{\infty}\,dt_k\int_{t_j\wedge
t}^{t_{j+1}\wedge t}\,dt^{\prime}
\int_E q(y_0,dy_1)\cdots\\
&&\qquad\quad\int
_E q(y_{j-1}, dy_j)\int
_E q\bigl(y_j,dy^{\prime}\bigr)
\\
&&\qquad\quad \int_E q\bigl(\pi_{[0,s-\varepsilon]}^E
\bigl(y^{\prime}\bigr),dy_{j+1}\bigr)\cdots\int
_E q(y_{k-1}, dy_k)I_{B}(t_1,y_1,
\ldots,t_k,y_k)I_{J_{y_j}}\bigl(y^{\prime}
\bigr)
\\
&&\qquad\quad{} \cdot\exp \Biggl\{-\sum_{l=0}^{j-1}q
\bigl(y_l,\bigl[y_l^s\bigr]^{c}
\bigr) (t_{l+1}-t_l)- \sum_{l=j+1}^{k-1}q
\bigl(y_l,\bigl[y_l^{s-\varepsilon}\bigr]^{c}
\bigr) (t_{l+1}-t_l)
\\
&&\hspace*{62pt}\qquad\quad{} -q\bigl(y_j,\bigl[y_j^s
\bigr]^{c}\bigr) \bigl(t^{\prime}-t_j\bigr)-q
\bigl(y^{\prime},\bigl[\bigl(y^{\prime}\bigr)^{s-\varepsilon}
\bigr]^{c}\bigr) \bigl(t_{j+1}-t^{\prime}\bigr) \Biggr
\}.
\end{eqnarray*}
Since $y_j\in E^a$, hence for $y^{\prime}\in J_{y_j}$ we have $\pi
_{[0,s-\varepsilon]}^E(y^{\prime})=\pi_{[0,s-\varepsilon]}^E(y_j)=y_j$ and $\int_E
q(y_j,dy^{\prime})I_{J_{y_j}}(y^{\prime})=\int_{s-\varepsilon}^{s}2^{-1}\rho
p(u)\,du$. Therefore,
\begin{eqnarray*}
&&P(K_{a,s,\varepsilon,j,\xi}\cap H\cap A_{k,a,j,\xi})\\
&&\qquad=\int_0^{\infty}\,dt_1\cdots\int
_0^{\infty}\,dt_k\int_E
q(y_0,dy_1) \cdots \int_E
q(y_{k-1}, dy_k)
\\
&&\qquad\quad I_{B}(t_1,y_1,\ldots,t_k,y_k)
\int_{t_j\wedge t}^{t_{j+1}\wedge
t}\,dt^{\prime}\int
_{s-\varepsilon}^{s}2^{-1}\rho p(u)\,du
\\
&&\qquad\quad{} \cdot\exp \Biggl\{-\sum_{l=0}^{j-1}q
\bigl(y_l,\bigl[y_l^s\bigr]^{c}
\bigr) (t_{l+1}-t_l) -\sum_{l=j+1}^{k-1}q
\bigl(y_l,\bigl[y_l^{s-\varepsilon}\bigr]^{c}
\bigr) (t_{l+1}-t_l)
\\[-2pt]
&&\hspace*{72pt}\qquad\quad{} -q\bigl(y_j,\bigl[y_j^s
\bigr]^{c}\bigr) \bigl(t^{\prime}-t_j\bigr)-q
\bigl(y_j,\bigl[y_j^{s-\varepsilon}\bigr]^{c}
\bigr) \bigl(t_{j+1}-t^{\prime}\bigr) \Biggr\}.
\end{eqnarray*}
For $ (t_1,y_1,\ldots,t_k,y_k)\in B$, we have $y_l\in E^a$ for all
$1\leq l \leq k$, therefore,\break
$q(y_l,[y_l^s]^{c})=q(y_l,[y_l^{s-\varepsilon}]^{c})+|\pi_a(y_l)|
\int_{s-\varepsilon}^{s}2^{-1}\rho p(u)\,du$. Thus,
\begin{eqnarray*}
&&P(K_{a,s,\varepsilon,j,\xi}\cap H\cap A_{k,a,j,\xi})\\[-2pt]
&&\qquad=\int_0^{\infty}\,dt_1\cdots\int
_0^{\infty}\,dt_k\int_E
q(y_0,dy_1)\cdots \int_E
q(y_{k-1}, dy_k)
\\[-2pt]
&&\qquad\quad I_{B}(t_1,y_1,\ldots,t_k,y_k)
\int_{s-\varepsilon}^{s}2^{-1}\rho p(u)\,du\int
_{t_j\wedge t}^{t_{j+1}\wedge t}\,dt^{\prime}
\\[-2pt]
&&\qquad\quad{} \cdot\exp \Biggl\{ -\sum_{l=0}^{j-1}
\biggl[q\bigl(y_l,\bigl[y_l^{s-\varepsilon}
\bigr]^{c}\bigr)+\bigl| \pi_a(y_l)\bigr|\int
_{s-\varepsilon}^{s}2^{-1}\rho p(u)\,du
\biggr](t_{l+1}-t_l)
\\[-2pt]
&&\hspace*{28pt}\qquad\quad{} -\sum_{l=j+1}^{k-1}q\bigl(y_l,
\bigl[y_l^{s-\varepsilon}\bigr]^{c}\bigr)
(t_{l+1}-t_l)-q\bigl(y_j,\bigl[y_j^{s-\varepsilon}
\bigr]^{c}\bigr) \bigl(t_{j+1}-t^{\prime}\bigr)
\\[-2pt]
&&\hspace*{46pt}\qquad\quad{} - \biggl[q\bigl(y_j,\bigl[y_j^{s-\varepsilon}
\bigr]^{c}\bigr)+\bigl|\pi_a(y_j)\bigr|\int
_{s-\varepsilon
}^{s}2^{-1}\rho p(u)\,du \biggr]
\bigl(t^{\prime}-t_j\bigr) \Biggr\}.
\end{eqnarray*}
Note that $q(y_l,[y_l^{s-\varepsilon}]^{c})=q(y_l,[y_l^a]^{c})+|\pi_a(y_l)|
\int_{a}^{s-\varepsilon}2^{-1}\rho p(u)\,du$. We define $\tilde{L}_a:=\sum_{l=0}^{k-1}(t_{l+1}-t_l)|\pi_a(y_l)|$, then
\begin{eqnarray*}
&&P(K_{a,s,\varepsilon,j,\xi}\cap H\cap A_{k,a,j,\xi})\\[-2pt]
&&\qquad=
\int_0^{\infty}\,dt_1\cdots\int
_0^{\infty}\,dt_k\int_E
q(y_0,dy_1)\cdots \int_E
q(y_{k-1}, dy_k)
\\[-2pt]
&&\qquad\quad I_{B}(t_1,y_1,\ldots,t_k,y_k)
\int_{s-\varepsilon}^{s}2^{-1}\rho p(u)\,du \int
_{t_j\wedge t}^{t_{j+1}\wedge t}\,dt^{\prime}
\\[-2pt]
&&\qquad\quad{} \cdot\exp \Biggl\{-\sum_{l=0}^{k-1}q
\bigl(y_l,\bigl[y_l^a\bigr]^{c}
\bigr) (t_{l+1}-t_l)-2^{-1}\rho
\tilde{L}_a\int_{a}^{s-\varepsilon}p(u)\,du
\\[-2pt]
&&\hspace*{26pt}\qquad\quad{} -\sum_{l=0}^{j-1}\bigl|\pi_a(y_l)\bigr|(t_{l+1}-t_l)
\int_{s-\varepsilon
}^{s}2^{-1}\rho p(u)\,du
\\[-2pt]
&&\hspace*{95pt}\qquad\quad{} -\bigl|\pi_a(y_j)\bigr|\bigl(t^{\prime}-t_j
\bigr)\int_{s-\varepsilon
}^{s}2^{-1}\rho p(u)\,du
\Biggr\}.
\end{eqnarray*}
Multiplying with $ (1-\exp\{-\tilde{L}_{a} \int_{s-\varepsilon
}^{s}2^{-1}\rho p(u)\,du\} ) /  (1-\exp\{-\tilde{L}_{a}\times\break \int_{s-\varepsilon}^{s}2^{-1}\rho p(u)\,du\} )$ at the right-hand side of
the above equality, we get
\begin{eqnarray*}
&&P(K_{a,s,\varepsilon,j,\xi}\cap H\cap A_{k,a,j,\xi})\\
&&\qquad=\int_0^{\infty}\,dt_1\cdots\int
_0^{\infty}\,dt_k\int_E
q(y_0,dy_1)\cdots \int_E
q(y_{k-1}, dy_k)I_{B}(t_1,y_1,
\ldots,t_k,y_k)
\\
&&\qquad\quad{} \cdot\exp\Biggl\{-\sum_{l=0}^{k-1}q
\bigl(y_l,\bigl[y_l^a\bigr]^{c}
\bigr) (t_{l+1}-t_l)\Biggr\} \tilde{C}(t_1,y_1,
\ldots,t_k,y_k)
\\
&&\qquad\quad{} \cdot \biggl[\exp\biggl(-2^{-1}\rho\tilde{L}_a\int
_{a}^{s-\varepsilon
}p(u)\,du\biggr)-\exp\biggl(-2^{-1}
\rho\tilde{L}_a\int_{a}^{s}p(u)\,du
\biggr) \biggr],
\end{eqnarray*}
where
\begin{eqnarray*}
&&\tilde{C}(t_1,y_1,\ldots,t_k,y_k)\\
&&\qquad:=
\biggl(1-\exp\biggl\{-\tilde{L}_{a} \int_{s-\varepsilon}^{s}2^{-1}
\rho p(u)\,du\biggr\} \biggr)^{-1} \int_{s-\varepsilon
}^{s}2^{-1}
\rho p(u)\,du
\\
&&\qquad\quad {}\cdot\int_{t_j\wedge t}^{t_{j+1}\wedge t}\exp \Biggl
\{- \Biggl[\bigl|
\pi_a(y_j)\bigr|\bigl(t^{\prime}-t_j
\bigr)+ \sum_{l=0}^{j-1}\bigl|\pi
_a(y_l)\bigr|(t_{l+1}-t_l) \Biggr]
\\
&&\hspace*{179pt}\qquad\quad{}\cdot\int_{s-\varepsilon}^{s}2^{-1}\rho p(u)\,du
\Biggr\}\,dt^{\prime}.
\end{eqnarray*}
For $\omega=(t_1,y_1,\ldots,t_k,y_k)\in B$, if we set $g_\omega:
=y_0I_{[0,t_1)}(t)+\sum_{l=1}^{k}y_lI_{[t_l,t_{l+1})}(t)$, with the
convention that $t_{k+1}=\infty$, then $g_\omega\in G^a$. One can check
that $|\mathcal{T}_{a}(\tau_l^a)|(g_\omega)=|\pi_a(y_l)|,   L_a(g_\omega
)=\sum_{l=0}^{k-1}|\pi_a(y_{l})|(t_{l+1}-t_{l})=\tilde{L}_a$ and
$C(a,j,\break s, \varepsilon)(g_\omega)=\tilde{C}(t_1,y_1,\ldots,t_k,y_k)$.
Therefore, applying Proposition~\ref{distribution2} we obtain
\begin{eqnarray*}
&&P(K_{a,s,\varepsilon,j,\xi}\cap H\cap A_{k,a,j,\xi})
\\
&&\qquad=\int_{H}C(a,j,s,\varepsilon) \biggl[\exp
\biggl(-2^{-1}\rho{L}_a\int_{a}^{s-\varepsilon}p(u)\,du
\biggr)\\
&&\hspace*{106pt}{}-\exp\biggl(-2^{-1}\rho{L}_a\int_{a}^{s}p(u)\,du
\biggr) \biggr]A_{k,a,j,\xi} P(dg).
\end{eqnarray*}
Since $H\in\sigma(X^a)$ is arbitrary, hence what we have proved
implies that
\begin{eqnarray*}
&&P\bigl(K_{a,s,\varepsilon,j,\xi}\cap H\cap A_{k,a,j,\xi
} |X^a\bigr)
\\
&&\qquad= C(a,j,s,\varepsilon) \biggl[\exp\biggl(-2^{-1}\rho{L}_a
\int_{a}^{s-\varepsilon
}p(u)\,du\biggr)
 \\
 &&\hspace*{91pt}{}-\exp\biggl(-2^{-1}\rho{L}_a\int_{a}^{s}p(u)\,du
\biggr) \biggr]A_{k,a,j,\xi}
\\
&&\qquad=C(a,j,s,\varepsilon)P\bigl(\{s-\varepsilon\leq S_{i+1}< s\}
|X^a\bigr)A_{k,a,j,\xi},
\end{eqnarray*}
where the last line is due to Lemma~\ref{Si2}.
The proof is completed by noticing that $A_{k,a,j,\xi}\in\sigma(X^a)$.
\end{pf}
For $n\geq1$, we define $\sigma_n=2^{-n}([2^nS_{i+1}]+1)$. Then each
$\sigma_n$ is a countably valued $\{\mathcal{F}^s\}$ stopping time and
$\sigma(S_{i+1})=\bigvee_n\sigma(\sigma_n)$. If we replace $s$ by $\sigma
_n$ and $\varepsilon$ by $2^{-n}$ in (\ref{K1}), we get another subset
$\tilde{K}_{a,\sigma_n,j,\xi}$ from the expression of $K_{a,s,\varepsilon,j,\xi}$ as follows:
%
%e4.9 #&#
\begin{eqnarray}\qquad
&&\tilde{K}_{a,\sigma_n,j,\xi}=\bigl\{\tau_{j+1}^{\sigma
_n}\leq t,
\tau_{j}^a<\tau_{j+1}^{\sigma_n}<
\tau_{j+1}^{a}, \exists r\in\bigl(\sigma_n-2^{-n},
\sigma_n\bigr)
\nonumber
\\[-8pt]
\\[-8pt]
\nonumber
&&\hspace*{56pt}\mbox{and } f\in X^{\sigma_n}\bigl(\tau _{j+1}^{\sigma_n}
\bigr), \mbox{ s.t. }f(r)=\xi, f(u)=\varnothing \mbox{ for all } u<r\bigr\}.
\nonumber
\end{eqnarray}

%le5 #&#
\begin{lemma}\label{distn}
Let the notation be the same as the above lemma. For any $n\geq1$ and
$t>0$, we have
\begin{eqnarray*}
&&P\bigl(\tilde{K}_{a,\sigma_n,j,\xi} |X^a, \sigma_n\bigr)
A_{k,a,j,\xi}I_{\{a< \sigma_n-{1}/{2^n}\}}
\\
&&\qquad=C\biggl(a,j,\sigma_n,\frac{1}{2^n}\biggr)A_{k,a,j,\xi}
I_{\{a< \sigma_n-{1}/{2^n}\}}.
\end{eqnarray*}
\end{lemma}
\begin{pf}
By the definition of $\sigma_n$, we have $\sigma_n=\sum_{m\geq1}\frac
{m}{2^n}I\{{m}/{2^n}-{1}/{2^n}\leq  S_{i+1}< {m}/{2^n}\}
$, hence $\tilde{K}_{a,\sigma_n,j,\xi}=\bigcup_{m\geq1}K_{a,
{m}/{2^n},{1}/{2^n},j,\xi}$. Then by Lemma~\ref{distZ1} one can
check directly that
\begin{eqnarray*}
&&P\bigl(\tilde{K}_{a,\sigma_n,j,\xi} |X^a, \sigma_n\bigr)
A_{k,a,j,\xi}I_{\{a< \sigma_n-{1}/{2^n}\}}
\\
&&\qquad=\sum_{m\geq1}P \bigl(K_{a,{m}/{2^n},{1}/{2^n},j,\xi}|X^a
\bigr) \biggl(P\biggl(\biggl\{\frac{m}{2^n}-\frac{1}{2^n}\leq
S_{i+1}< \frac{m}{2^n}\biggr\} \Big|X^a\biggr)
\biggr)^{-1}
\\
&&\qquad\quad{} \cdot A_{k,a,j,\xi}I_{\{a< {m}/{2^n}-{1}/{2^n}\}}
I_{\{{m}/{2^n}-{1}/{2^n}\leq S_{i+1}< {m}/{2^n}\}
}
\\
&&\qquad=\sum_{m\geq1}C\biggl(a,j,\frac{m}{2^n},
\frac{1}{2^n}\biggr)A_{k,a,j,\xi}I_{\{a<
{m}/{2^n}-{1}/{2^n}\}}I_{\{\sigma_n={m}/{2^n}\}}
\\
&&\qquad=C\biggl(a,j,\sigma_n,\frac{1}{2^n}\biggr)A_{k,a,j,\xi}
I_{\{a< \sigma_n-{1}/{2^n}\}}.
\end{eqnarray*}
\upqed\end{pf}

\begin{pf*}{Proof of Theorem~\ref{distT0}}
For $m\geq1$, we define $\alpha_m=2^{-m}([2^mS_{i}]+1)$. If we replace
$a$ by $\alpha_m$ in $\tilde{K}_{a,\sigma_n,j,\xi}$, we get another
subset $\tilde{K}_{\alpha_m,\sigma_n,j,\xi}$.
Similar to the proof of Lemma~\ref{distn}, we can show that
\begin{eqnarray*}
&&P\bigl(\tilde{K}_{\alpha_m,\sigma_n,j,\xi} |X^{\alpha_m}, \sigma_n\bigr)
A_{k,\alpha_m,j,\xi}I_{\{\alpha_m< \sigma_n-{1}/{2^n}\}}
\\
&&\qquad=C\biggl(\alpha_m,j,\sigma_n,\frac{1}{2^n}
\biggr)A_{k,\alpha_m,j,\xi}I_{\{\alpha
_m< \sigma_n-{1}/{2^n}\}}.
\end{eqnarray*}
Since $\lim_{m\rightarrow\infty}\alpha_m=S_i$ and $\sigma(S_i)=\bigvee
_m\sigma(\alpha_m)$, by the convergence property of conditional
expectations (cf., e.g., \cite{Yan92}, Theorem~2.21), we get
%
%e4.10 #&#
\begin{eqnarray}
\label{kj} &&P\bigl(\tilde{K}_{S_i,\sigma_n,j,\xi} |X^{S_i},
\sigma_n\bigr) A_{k,S_i,j,\xi}I_{\{S_i< \sigma_n-{1}/{2^n}\}}
\nonumber
\\[-8pt]
\\[-8pt]
\nonumber
&&\qquad=C\biggl(S_i,j,\sigma_n,\frac{1}{2^n}
\biggr)A_{k,S_i,j,\xi}I_{\{S_i< \sigma_n-{1}/{2^n}\}}.
\end{eqnarray}
Note that by (\ref{Cajs}) we have
\[
\lim_{n\rightarrow\infty}C\biggl(S_i,j,\sigma_n,
\frac{1}{2^n}\biggr)= \frac
{1}{L_{S_i}} \int_{\tau_j^{S_i}\wedge t}^{\tau_{j+1}^{S_i}\wedge t}\,du.
\]
By the definition of $A_{k,S_i,j,\xi}$, we have
\[
A_{k,S_i,j,\xi}=\bigl\{ \gamma^{S_i}=k, \xi\in\mathcal{T}_{S_i}
\bigl(\tau _{j}^{S_i}\bigr)\bigr\}.
\]
Noticing that $\tilde{K}_{S_i,\sigma_n,j,\xi}\rightarrow\{
T_0^{i+1}\leq t, \xi^{i+1}=\xi,T^{i+1}_0\in(\tau^{S_i}_j,\tau
^{S_i}_{j+1})\}$ a.s. as $n\rightarrow\infty$, letting $n\rightarrow
\infty$ in the both sides of (\ref{kj}) and employing again the
convergence property of conditional expectations, we get
\begin{eqnarray*}
&&P\bigl\{T_0^{i+1}\leq t, \xi^{i+1}=
\xi,T^{i+1}_0\in\bigl(\tau^{S_i}_j,\tau
^{S_i}_{j+1}\bigr)|X^{S_i},S_{i+1}\bigr
\}I_{\{ \gamma^{S_i}=k,  \xi\in\mathcal
{T}_{S_i}(\tau_{j}^{S_i})\}}
\\
&&\qquad=\frac{1}{L_{S_i}} \int_{\tau_j^{S_i}\wedge t}^{\tau_{j+1}^{S_i}\wedge
t}\,du\,
I_{\{ \gamma^{S_i}=k,  \xi\in\mathcal{T}_{S_i}(\tau_{j}^{S_i})\}}.
\end{eqnarray*}
Summing up the above equation for $k$ and $j$, and noticing that $\xi\in
\mathcal{T}_{S_i}(\tau_{j}^{S_i})$ if and only if $\xi\in\mathcal
{T}_{S_i}(u)$ for all $\tau_{j}^{S_i} \leq u < \tau_{j}^{S_{i+1}}$, we
complete the proof of the theorem.
\end{pf*}

%s4.5 #&#
\subsection{Proof of Proposition \texorpdfstring{\protect\ref{finitedim}}{16}}\label{prooffinitedim}

By Theorem~\ref{distZ^0}, we see that the distribution of $Z^0$
generated by step 1 coincides the distribution of $Z^0$ defined by (\ref
{Z^O}). It is apparent that the conditional distribution of $S_{i+1}$
generated by step 2 coincides with the one specified by Theorem~\ref
{distSi}, and
the conditional distribution of $T_0^{i+1}$ and $\xi^{i+1}$ generated
by step 3 is the same as those described by Theorem~\ref{distT0}. To
analyze the random
elements generated in steps 4 and~5, we define recursively for
$n\geq1$,
%
%e4.11 #&#
\begin{equation}
\label{T^in} T^i_n= \inf\bigl\{t>T^i_{n-1}
\dvtx Z^i(t)\neq Z^i\bigl({T^i_{n-1}}
\bigr)\bigr\}\quad\mbox{and}\quad\xi^i_n=Z^i
\bigl(T^i_n\bigr).
\end{equation}

In step 4, it is implicitly assumed that the ancestral material carried
on the new branch is $\xi_0^{i+1}$, which means that $\xi_j^{i+1}=\xi_0^{i+1}$.
It is not difficult to check that the distribution employed in step 4
coincides with the one developed in Theorem~\ref{distTxi1}. To analyze
step 5, we note that step 5 corresponds to the case that $\xi
_j^{i+1}\neq\xi_0^{i+1}$, and the EDGE is labeled with $i$ if and only
if $\pi_i(\xi_j^{i+1})\neq\varnothing$. When the EDGE is labeled with
$i$, the algorithm goes to step 6, the path
carrying the ancestral material $\xi_0^{i+1}$ is assumed to move along
the edges of $X^{S_{i}}$. In this case, $T_{j+1}^{i+1}$ must be the
jump time of $X^{S_i}$ at which the lineage carrying the ancestral
material $\pi_{[0,i]}(\xi_j^{i+1})$ meets its first change after $T_{j}^{i+1}$.
Thus, the conditional distribution of $(T^{i+1}_n, \xi_n^{i+1})$
coincides with the one described in Theorem~\ref{distTxi3}. In step 5,
when the EDGE is labeled with some $k$ less than $i$, then a potential
recombination event is considered in the algorithm. We point out that
$k$ is equal to $h(\vec{\xi_j})$ used in Theorem~\ref{distTxi2}, and
the upper node of the EDGE, denoted by $\mathcal{H}$ in the algorithm,
is the time point $\mathcal{H}$ used in Theorem~\ref{distTxi2}. Then
one can check that the distribution used in step 5 coincides with the
conditional distribution developed in Theorem~\ref{distTxi2}. To sum up
the above discussion, we find that all the distributions of the random
elements generated by the algorithm coincide with those developed in
Section~\ref{sequence}. Therefore, the finite dimensional
distribution of the random sequence
$\{(S_i, Z^i),i\geq0\}$ generated by the $\mathit{SC}$ algorithm is the same as
that developed in Section~\ref{sequence}.

%s4.6 #&#
\subsection{Proof of Theorem \texorpdfstring{\protect\ref{coincide}}{8}}\label{proofcoincide}

Below we use $\{(\tilde{S_i}, \tilde{Z^i}),i\geq0\}$ to denote the
random elements $\{(S_i, Z^i),i\geq0\}$ generated by the $\mathit{SC}$
algorithm, and reserve $\{(S_i, Z^i),i\geq0\}$ for those originally
defined on $(G,\mathcal{B}(G),P)$ as discussed in Section~\ref{sequence}. It is implicitly assumed that $\{(\tilde{S_i}, \tilde
{Z^i}),i\geq0\}$ are taken from some probability space other than the
space $G$. We denote by
%
%e4.12 #&#
\begin{equation}
\tilde{\Phi}_i:=\bigl(\tilde{S_1},\ldots,
\tilde{S_i};\tilde{Z^0},\tilde {Z^1},\ldots,
\tilde{Z^i}\bigr),
\end{equation}
and denote by $\tilde{P}_i$ the probability distribution of $\tilde{\Phi
}_i$ on its sample space $\Omega_i$ specified by (\ref{Omigai}). Here,
\[
\Omega_i:=[0,1]^{i} \times\mathcal{S}_{[0,\infty)}(
\mathcal{R}) \times \prod_{l=1}^{i}
\mathcal{S}_{[0,\infty)}\bigl(\mathcal{P}^{l+1}\bigr)
\]
is also the sample space of $\Phi_i:=(S_1,\ldots,S_i;Z^0,Z^1,\ldots,Z^i)$
specified by (\ref{Phii}). Denote by $P_i$ the probability distribution
of $\Phi_i$. By Proposition~\ref{finitedim}, we have $\tilde{P}_i=P_i$.
Bellow we use the notation employed in the proof of Proposition~\ref{fil}.
Then $\tilde{P}_i=P_i$ implies in particular that $\tilde{P}_i(\mathcal
{H}_i)=1$. Each $\omega\in\mathcal{H}_i$ constitutes a
part graph $g \in X^{S_i}(G):=G_i$, which is described by the map
$\Upsilon_i: =(\Phi_i|_{G_i})^{-1}$. Since $\Upsilon_i\dvtx \mathcal{H}_i
\mapsto G_i$ is a one to one Borel map, hence $\tilde{P}_i$ induces a
probability measure $\tilde{P}_i^*=\tilde{P}_i\circ(\Upsilon_i)^{-1}$
on $G_i$. Similarly, $P_i$ induces a
probability measure $P_i^*$ on $G_i$ and we have $\tilde{P}_i^*=P_i^*$.
By Lemma~\ref{F^Si}, we have $\sigma(X^{S_i})=\mathcal{F}^{S_i}$. For
notational convenience below, we write $\mathcal{F}_i$ for $\mathcal
{F}^{S_i}$. Through the mapping $X^{S_i}\dvtx G_i \mapsto G$, the
probability measure $P_i^*$ determines a probability measure
$P^*_{|\mathcal{F}_i}$ on $(G,\mathcal{F}_i)$ by setting $P^*_{|\mathcal
{F}_i}((X^{S_i})^{-1}(B))=P^*_i(B)$ for all $B\in\mathcal{B}(G_i)$.
Noticing that $X^{S_i}(g)=\Upsilon_i(\Phi_i(g))$, one can check that
$P^*_{|\mathcal{F}_i}=P_{|\mathcal{F}_i}$ where $P_{|\mathcal{F}_i}$ is
the restriction of $P$ on $\mathcal{F}_i$. Similarly, $\tilde{P}_i^*$
determines a probability measure $\tilde{P}^*_{|\mathcal{F}_i}$ on
$(G,\mathcal{F}_i)$.
On the other hand, let $\tilde{P}$ be the probability distribution on
$(G,\mathcal{B}(G))$
generated by the algorithm $\mathit{SC}$, then we must have $\tilde{P}\circ
(X^{S_i})^{-1}=\tilde{P}_i^*$ on $G_i$. Then it is clear that $\tilde
{P}^*_{|\mathcal{F}_i}=\tilde{P}_{|\mathcal{F}_i}$.
Therefore, we get $\tilde{P}_{|\mathcal{F}_i}=\tilde{P}^*_{|\mathcal
{F}_i}=P^*_{|\mathcal{F}_i}
=P_{|\mathcal{F}_i}$. The proof of the theorem is completed by noticing that
$\mathcal{B}(G)=\bigvee_{i\geq1}\mathcal{F}_i$.

\section*{Akonowledgments}
We are indebted to Shuhua Xu, Ying Zhou, Linfeng Li and Yuting Liu for
allowing us to use some material from our joint work \cite{wang}. We
are grateful to De-Xin Zhang and Wei-Wei Zhai for offering their
stimulating idea and for their encouragement. We thank
Renming Song for his very helpful comments and suggestions. We thank
the Editor, Associate Editor and the anonymous referees for their
valuable comments and suggestions which improved the presentation of
this paper.

\begin{supplement}[id=suppA]
\stitle{Supplement to ``Markov jump processes in modeling coalescent
with recombination''}
\slink[doi]{10.1214/14-AOS1227SUPP} %[doi,text={...}] - jei reikia
%suskaldyti doi
\sdatatype{.pdf}
\sfilename{aos1227\_supp.pdf}
\sdescription{The supplementary file is divided into two Appendixes.
Appendix A contains the proofs of Propositions~\ref{prop1}--\ref{3Gneighborhood} and Propositions
\ref{Zjcon}--\ref{measG}.
Appendix B is devoted to the calculation of the conditional
distribution $P(T_{j+1}^{i+1}\in
B, \xi_{j+1}^{i+1}=\vec{\xi} | X^{S_i},S_{i+1}, T_0^{i+1},
\xi^{i+1},\ldots,\break  T_j^{i+1},\xi_j^{i+1})$. In particular, the proofs of
Theorems \ref{distTxi1}, \ref{distTxi3} and~\ref{distTxi2}
are presented, respectively, in the proofs of Theorems B.10, B.11
and~B.12 in Appendix B.}
\end{supplement}

% imsref loaded by akundreckaite, 2014-05-27 08:23:06
% imsref loaded by akundreckaite, 2014-05-27 08:41:09


\printaddresses
\end{document}